\documentclass{conm-p-l}
\usepackage{amssymb,latexsym}
\usepackage{epic,eepic}
\usepackage{graphicx}
\usepackage[all]{xy}

\numberwithin{equation}{section}

\newtheorem{theorem}{Theorem}[section]
\newtheorem{proposition}[theorem]{Proposition}
\newtheorem{conjecture}[theorem]{Conjecture}
\newtheorem{corollary}[theorem]{Corollary}

\theoremstyle{definition}

\newtheorem{example}[theorem]{Example}
\newtheorem{definition}[theorem]{Definition}

\begin{document}
\bibliographystyle{amsalpha}

\title{Tropicalization method in cluster algebras}

\author{Tomoki Nakanishi}
\address{\noindent Graduate School of Mathematics, Nagoya University, 
Chikusa-ku, Nagoya,
464-8604, Japan}
\email{nakanisi@math.nagoya-u.ac.jp}

\begin{abstract}
This is a brief survey on
the recently developing tropicalization method
in cluster algebras
and its applications to the periodicities of Y-systems
 and the associated dilogarithm identities. 
\end{abstract}



\maketitle

\section{Introduction}

The nature of tropicalization is built into
cluster algebras from the beginning by the series of
fundamental works by Fomin, Zelevinsky, and Berenstein
\cite{Fomin02,Fomin03a,Fomin03b,Berenstein05,
Berenstein05b,Fomin07}.

Recall that there are two kinds of variables in cluster algebras,
namely, {\em cluster variables\/} and {\em coefficients}.
The cluster variables are generators of a cluster algebra itself,
and  they are the main characters in many applications of cluster algebras.
Meanwhile, the coefficients are a little bit
in a subsidiary position at first sight,
but they turn out to be equally important.
The coefficients take values in any choice of semifield.
In particular, one can choose it to be a {\em tropical\/} semifield.
Such a cluster algebra is said to be {\em of geometric type}
in \cite{Fomin02},
and this is where the tropicalization gets into business.

These `tropical' coefficients were introduced and
studied by several reasons:
\begin{itemize}
\item[i)] Tropical coefficients are also regarded as {\em frozen\/} cluster
variables (together with {\em extended\/} exchange matrices), and they naturally arise in various contexts in applications
\cite{Fomin02}.

\item[ii)] They appeared in the study
of the periodicity conjecture of Y-systems,
and  the dynamics of these tropical
coefficients are governed by the underlying Weyl group
and root system \cite{Fomin03b}.
This result  is further used in the classification of the
cluster algebras of finite type \cite{Fomin03a}.

\item[iii)] A special case of tropical coefficients,
called the {\em principal coefficients\/} play the central role 
of the structure theory of cluster variables and coefficients
\cite{Fomin07}.

\end{itemize}

Tropicalization also appeared (a little implicitly)
in the formulation of cluster algebras based on algebraic tori
by Fock and Goncharov
\cite{Fock03,Fock07,Fock07b}.
The importance of this approach becomes more manifest in
the analogous formulation of {\em quantum cluster algebras}
\cite{Fock03,Berenstein05b,Fock07,Fock07b,Tran09}.

More recently, the tropicalization method became more powerful
 with the help of the
{\em categorification\/} of cluster algebras by
{\em generalized cluster categories} developed by
\cite{Caldero06,Buan06,Dehy08,Fu10,Palu08,Amiot09,Keller08,Keller09,Plamondon10a,Plamondon10b},
etc.
In particular,
 we have two remarkable applications of the tropicalization method,
namely,
the periodicity of seeds
\cite{Keller08,Keller10,Inoue10a,Inoue10b,Inoue10c,Nakanishi10b}
and the associated dilogarithm identities
\cite{Chapoton05,Nakanishi09,Inoue10b,Inoue10c,Nakanishi10b,
Nakanishi10c}.
As corollaries, they solved two long standing conjectures
which arose in the study of integrable models in 1990's,
namely,
the {\em periodicities of Y-systems} and the {\em dilogarithm identities
 in conformal field theory}.
Furthermore,
 these dilogarithm identities have  their quantum counterparts
\cite{Nagao10,Keller11,Nagao11b,Nagao11,Kashaev11},
and the tropicalization plays a crucial role also in this context.

In this note we give a brief overview of these developments,
focusing on the role of tropicalization.
Essentially all the results are collected from the literature;
however, some part of the presentation (e.g., Section 
\ref{subsec:ytropx}) is new.

We thank Gen Kuroki, Takao Yamazaki and Andrei Zelevinsky for useful comments.

\section{Cluster algebra with coefficients}

We recall some basic facts on cluster algebras with coefficients in
\cite{Fomin02,Fomin03a,Fomin07} with a little change
of notation and terminology therein.
Throughout this section we fix a finite index set $I$,
say, $I=\{1,\dots,n\}$.

\subsection{Semifield}

\begin{definition}
We say that $\mathbb{P}$ is a {\em semifield\/}
if $\mathbb{P}$ is a multiplicative abelian group
endowed with the addition $\oplus$, which
is commutative, associative, and distributive,
i.e., $(a\oplus b) c = ac \oplus bc$.
\end{definition}

The following two examples of semifields are important
in this note.

\begin{example}
(a) {\em Universal semifield\/ $\mathbb{P}_{\mathrm{univ}}(y)$.}
For an $I$-tuple of variables $y=(y_i)_{i\in I}$,
it is the semifield consisting of all the rational functions of 
$y$ over $\mathbb{Q}$ with subtraction-free expressions.

\par
(b) {\em Tropical semifield\/ $\mathbb{P}_{\mathrm{trop}}(y)$.}
For an $I$-tuple of variables $y=(y_i)_{i\in I}$,
it is the multiplicative free abelian group
(i.e., the group of all the Laurent monomials of $y$ with coefficient 1)
endowed with the following addition $\oplus$:
\begin{align}
\label{eq:tropo}
  \prod_{i\in I} y_i ^{a_i}
\oplus
\prod_{i\in I} y_i ^{b_i}
:=
\prod_{i\in I} y_i ^{\min (a_i,b_i)}.
\end{align}

Later we use the following natural homomorphism of semifields
\begin{align}
\label{eq:hom}
\pi_{\mathrm{trop}}:
\mathbb{P}_{\mathrm{univ}}(y)
\rightarrow
\mathbb{P}_{\mathrm{trop}}(y),
\end{align}
which sends $y_i \mapsto y_i$ and $c \mapsto 1$ ($c\in \mathbb{Q}_{> 0}$).
\end{example}

\subsection{Cluster algebra with coefficients}

Let $(B,x,y)$ be a triplet with the following data:
\begin{itemize}
\item
a {\em skew-symmetrizable\/} integer matrix
$B=(b_{ij})_{i,j\in I}$;
namely there is a positive integer
diagonal matrix $D=\mathrm{diag}(d_i)_{i\in I}$
such that $DB$ is skew-symmetric, i.e.,
$(DB)^T=-DB$, where $T$ is the transposition,
\item an $I$-tuple $x=(x_i)_{i\in I}$ of formal variables,
\item an $I$-tuple $y=(y_i)_{i\in I}$ of formal variables.
\end{itemize}
We fix $D$ such that
{\em $d_1,\dots,d_n$ are coprime\/} throughout.
Thus, $D=I$ when $B$ itself is skew-symmetric.
Below we consider the {\em 
cluster algebra with coefficients $\mathcal{A}(B,x,y)$},
whose initial seed is $(B,x,y)$ and
coefficients semifield is $\mathbb{P}_{\mathrm{univ}}(y)$.
Let us recall its definition  to fix the convention.

Let $(B',x',y')$ be a triplet with the following data:
\begin{itemize}
\item
a {\em skew-symmetrizable\/} integer matrix
$B'=(b'_{ij})_{i,j\in I}$;
such that $DB$ is skew-symmetric for the above $D$,

\item an $I$-tuple $x'=(x'_i)_{i\in I}$ 
with $x'_i \in \tilde{\mathbb{Q}}(x)$,
where $\tilde{\mathbb{Z}} :=\mathbb{Z}\mathbb{P}_{\mathrm{univ}}(y)$
is the group ring of $\mathbb{P}_{\mathrm{univ}}(y)$,
and $\tilde{\mathbb{Q}} :=\mathbb{Q}\mathbb{P}_{\mathrm{univ}}(y)$
is the fractional field of $\tilde{\mathbb{Z}}$.

\item an $I$-tuple $y'=(y'_i)_{i\in I}$ with
$y'_i \in \mathbb{P}_{\mathrm{univ}}(y)$.
\end{itemize}
Then, for $(B',x',y')$ and $k\in I$,
another triplet of the same kind
$(B'',x'',y'')$,
called the {\em mutation of $(B',x',y')$ at $k$}
and denoted by $\mu_k(B',x',y')$,
is defined by
the following rule:
\begin{align}
\label{eq:bmut}
b''_{ij}&=
\begin{cases}
-b'_{ij}& \mbox{$i=k$ or $j=k$}\\
b'_{ij}+[- b'_{ik}]_+ b'_{kj} + b'_{ik}[ b'_{kj}]_+
& \mbox{$i,j\neq k$,}\\
\end{cases}
\\
\label{eq:ymut}
y''_i&=
\begin{cases}
y'_k{}^{-1} & i= k\\
\displaystyle
y'_i \frac{(1\oplus y'_k)^{[-b'_{ki}]_+}}
{(1\oplus y'_k{}^{-1})^{[b'_{ki}]_+}}
& i\neq k.\\
\end{cases}
\\
\label{eq:xmut}
x''_i&=
\begin{cases}
\displaystyle
{x'_k}^{-1}\left(
\frac{1}{1\oplus y'_k{}^{-1}}
\prod_{j\in I} x'_j{}^{[b'_{jk}]_+}
+ 
\frac{1}{1\oplus y'_k}
\prod_{j\in I} x'_j{}^{[-b_{jk}]_+}
\right)
& i=k\\
x'_i & i\neq k.\\
\end{cases}
\end{align}
Here, for any integer $x$, we set $[x]_+:=\max(x,0)$.
The relations \eqref{eq:ymut} and \eqref{eq:xmut}
are called the {\em exchange relations}.
The involution property $\mu_k^2=\mathrm{id}$ holds.

Repeat mutations from the initial seed $(B,x,y)$
and collect all the obtained triplets $(B',x',y')$,
which are called {\em seeds}.
For each seed $(B',x',y')$, $B'$ is called an {\em exchange matrix},
$x'$ and $x'_i$ ($i\in I$) are called a {\em cluster\/} and
 a {\em cluster variable},
and $y'$ and  $y'_i$ ($i\in I$) are called a {\em coefficient tuple}
and a {\em coefficient}.
The {\em cluster algebra 
$\mathcal{A}(B,x,y)$ with coefficients in
$\mathbb{P}_{\mathrm{univ}}(y)$} is
the $\tilde{\mathbb{Z}}$-subalgebra of $\tilde{\mathbb{Q}}(x)$
generated by all the cluster variables.

In this note we casually call $x'_i$ and $y'_i$
the {\em $x$-variables\/} and the {\em $y$-variables}.
It is a little unfortunate and inconvenient that this conflicts with the
notation and terminology
(the $a$-variables and the $x$-variables)  used by another standard
and important references
by Fock and Goncharov \cite{Fock03,Fock07,Fock07b}.

For each seed $(B',x',y')$, define $\hat{y}'_i$ ($i\in I$) by
\begin{align}
\label{eq:yhat}
\hat{y}'_i:=y'_i\prod_{j\in I} x'_j{}^{b'_{ji}}.
\end{align}
Then, $\hat{y}'_i$ satisfies the same exchange relation 
\eqref{eq:ymut} as $y'_i$ \cite[Prop.3.9]{Fomin07};
 namely,
\begin{align}
\label{eq:yhatmut}
\hat{y}''_i=
\begin{cases}
\hat{y}'_k{}^{-1} & i= k\\
\displaystyle
\hat{y}'_i \frac{(1+ \hat{y}'_k)^{[-b'_{ki}]_+}}
{(1+ \hat{y}'_k{}^{-1})^{[b'_{ki}]_+}}
& i\neq k.\\
\end{cases}
\end{align}

\subsection{$\varepsilon$-expressions for mutations}

The matrix mutation \eqref{eq:bmut} is also written as
\begin{align}
\label{eq:bmut2}
b''_{ij}=
\begin{cases}
-b'_{ij}& \mbox{$i=k$ or $j=k$}\\
b'_{ij}+[- \varepsilon b'_{ik}]_+ b'_{kj} + b'_{ik}[ \varepsilon b'_{kj}]_+
& \mbox{$i,j\neq k$,}\\
\end{cases}
\end{align}
where $\varepsilon \in \{1, -1\}$ and it is {\em independent of the choice
of $\varepsilon$} \cite[Eq.(3.1)]{Berenstein05}. If $\varepsilon =1$, then
it is the same as \eqref{eq:bmut}.
In the same spirit, the exchange relations \eqref{eq:ymut}
and \eqref{eq:xmut}
 have the expression
(independent of the choice of $\varepsilon=\pm1$) as
\begin{align}
\label{eq:ymut2}
y''_i&=
\begin{cases}
y'_k{}^{-1} & i= k\\
\displaystyle
y'_i y'_k{}^{[\varepsilon b'_{ki}]_+}(1\oplus y'_k{}^{\varepsilon})^{-b'_{ki}}
& i\neq k,\\
\end{cases}\\
\label{eq:xmut2}
x''_i&=
\begin{cases}
\displaystyle
x'_k{}^{-1}
\left(
\prod_{j\in I} x'_j{}^{[-\varepsilon b'_{jk}]_+}
\right)
\frac{1+\hat{y}'_k{}^{\varepsilon}}{1\oplus y'_k{}^{\varepsilon}}
& i=k\\
x'_i & i\neq k,\\
\end{cases}
\end{align}
where $\hat{y}'_i$ is defined in \eqref{eq:yhat}.
Let us call \eqref{eq:bmut2}--\eqref{eq:xmut2} the
{\em $\varepsilon$-expressions\/} for their counterparts.

At first sight it seems useless to have such a redundant choice of
the sign $\varepsilon$.
However, it turns out that  there is actually a
{\em canonical choice of $\varepsilon$\/} at each mutation
{\em in view of tropicalization}.
We will see more about it below.

\subsection{Example}

It is standard to identify a  skew-symmetric matrix $B$ with a quiver
$Q$ (without loop and 2-cycle)
such that $I$ is the set of the vertices of $Q$, and 
$t$ arrows are drawn from the vertex $i$ to $j$ if
$b_{ij}=t >0$.

\begin{example}
\label{ex:a2}
Consider the cluster algebra $\mathcal{A}(B,x,y)$
whose initial exchange matrix $B$ and
the corresponding quiver $Q$
are given by
\begin{align}
B=
\begin{pmatrix}
0 & -1 \\
1 & 0\\
\end{pmatrix},
\quad
Q =\ 
\raisebox{-10pt}{
\begin{picture}(40,5)(0,-15)
%
%
\put(0,0){\circle{6}}
\put(40,0){\circle{6}}
\put(37,0){\vector(-1,0){34}}
\put(-3,-15){$1$}
\put(37,-15){$2$}
\end{picture}
}
\, .
\end{align}
This is the cluster algebra of type $A_2$ in
the classification of \cite{Fomin02}.
In particular, it is of finite type,
namely, there are only finitely many seeds.
Set $\Sigma(0)=(Q(0),x(0),y(0))$ to be the initial seed $(Q,x,y)$,
and consider the seeds
$\Sigma(t)=(Q(t), x(t), y(t))$ ($t=1,\dots,5$)
by the following sequence of the alternative mutations
of  $\mu_1$ and  $\mu_2$:
\begin{align}
\label{eq:seedmutseq3}
&\Sigma(0)
\
\mathop{\leftrightarrow}^{\mu_{1}}
\
\Sigma(1)
\
\mathop{\leftrightarrow}^{\mu_{2}}
\
\Sigma(2)
\
\mathop{\leftrightarrow}^{\mu_{1}}
\
\Sigma(3)
\
\mathop{\leftrightarrow}^{\mu_{2}}
\
\Sigma(4)
\
\mathop{\leftrightarrow}^{\mu_{1}}
\
\Sigma(5).
\end{align}
According to \eqref{eq:yhat}, we set
\begin{align}
\hat{y}_1 = y_1 x_2,
\quad \hat{y}_2=y_2 x_1^{-1}.
\end{align}
Then, using the exchange relations \eqref{eq:bmut2}--\eqref{eq:xmut2},
we obtain the following explict form of seeds:
\begin{alignat*}{3}
Q(0)\quad
&
\begin{picture}(40,15)(0,-5)
%
%
\put(0,0){\circle{12}}
\put(0,0){\circle{6}}
\put(40,0){\circle{6}}
\put(37,0){\vector(-1,0){34}}
\put(-3,-15){$1$}
\put(37,-15){$2$}
\end{picture}
&\quad&
\begin{cases}
x_1(0)=x_1\\
x_2(0)=x_2,\\
\end{cases}
&&
\begin{cases}
y_1(0)=y_1\\
y_2(0)=y_2,\\
\end{cases}
\notag
\\
Q(1)\quad
&
\begin{picture}(40,15)(0,-5)
%
%
\put(40,0){\circle{12}}
\put(0,0){\circle{6}}
\put(40,0){\circle{6}}
\put(3,0){\vector(1,0){34}}
\put(-3,-15){$1$}
\put(37,-15){$2$}
\end{picture}
&&
\begin{cases}
\displaystyle
x_1(1)=x_1^{-1}\frac{1+\hat{y}_1}{1\oplus y_1}\\
x_2(1)=x_2,\\
\end{cases}
&&
\begin{cases}
y_1(1)=y_1^{-1}\\
y_2(1)=y_2(1\oplus y_1),\\
\end{cases}
\\
Q(2)\quad
&
\begin{picture}(40,15)(0,-5)
%
%
\put(0,0){\circle{12}}
\put(0,0){\circle{6}}
\put(40,0){\circle{6}}
\put(37,0){\vector(-1,0){34}}
\put(-3,-15){$1$}
\put(37,-15){$2$}
\end{picture}
&&
\begin{cases}
\displaystyle
x_1(2)=x_1^{-1}\frac{1+ \hat{y}_1}{1\oplus y_1}\\
\displaystyle
x_2(2)=x_2^{-1}
\frac{1 + \hat{y}_2 +\hat{y}_1\hat{y}_2}{1\oplus y_2\oplus y_1y_2},\\
\end{cases}
&&
\begin{cases}
y_1(2)=y_1^{-1}(1\oplus y_2 \oplus y_1y_2)\\
y_2(2)=y_2^{-1}(1\oplus y_1)^{-1},\\
\end{cases}
\\
Q(3)\quad
&
\begin{picture}(40,15)(0,-5)
%
%
\put(40,0){\circle{12}}
\put(0,0){\circle{6}}
\put(40,0){\circle{6}}
\put(3,0){\vector(1,0){34}}
\put(-3,-15){$1$}
\put(37,-15){$2$}
\end{picture}
&&
\begin{cases}
\displaystyle
x_1(3)=x_1x_2^{-1}\frac{1+\hat{y}_2}{1\oplus y_2}
\\
\displaystyle
x_2(3)=
x_2^{-1}
\frac{1 + \hat{y}_2 +\hat{y}_1\hat{y}_2}{1\oplus y_2\oplus y_1y_2},\\
\end{cases}
&&
\begin{cases}
y_1(3)=y_1(1\oplus y_2 \oplus y_1y_2)^{-1}\\
y_2(3)=y_1^{-1}y_2^{-1}(1\oplus y_2),\\
\end{cases}
\\
Q(4)\quad
&
\begin{picture}(40,15)(0,-5)
%
%
\put(0,0){\circle{12}}
\put(0,0){\circle{6}}
\put(40,0){\circle{6}}
\put(37,0){\vector(-1,0){34}}
\put(-3,-15){$1$}
\put(37,-15){$2$}
\end{picture}
&&
\begin{cases}
\displaystyle
x_1(4)=x_1x_2^{-1}\frac{1+\hat{y}_2}{1\oplus y_2}\\
x_2(4)=x_1,\\
\end{cases}
&&
\begin{cases}
y_1(4)=y_2^{-1}\\
y_2(4)=y_1y_2(1\oplus y_2)^{-1},\\
\end{cases}
\\
Q(5)\quad
&
\begin{picture}(40,15)(0,-5)
%
%
\put(40,0){\circle{12}}
\put(0,0){\circle{6}}
\put(40,0){\circle{6}}
\put(3,0){\vector(1,0){34}}
\put(-3,-15){$1$}
\put(37,-15){$2$}
\end{picture}
&&
\begin{cases}
x_1(5)=x_2\\
x_2(5)=x_1.\\
\end{cases}
&&
\begin{cases}
y_1(5)=y_2\\
y_2(5)=y_1.\\
\end{cases}
\end{alignat*}
Here, the encircled vertices in quivers
are the {\em forward mutation points\/}
in the sequence \eqref{eq:seedmutseq3}.
We observe the (half) periodicity of mutations of seeds.
\end{example}

\subsection{$C$-matrices, $G$-matrices, and $F$-polynomials}

Let us present the  `separation formulas'
 due to
\cite{Fomin07}
clarifying the fundamental structure of $x$- and $y$-variables.

\begin{theorem}[{Separation formulas \cite[Prop.3.13 \& Cor.6.3]{Fomin07}}]
\label{thm:Fyx}
For each seed $(B',x',y')$ of any cluster algebra $\mathcal{A}(B,x,y)$,
there exist polynomials $F'_i(y)$ $(i \in I)$ of $y=(y_i)_{i\in I}$,
an integer matrix $C'=(c'_{ij})_{i,j\in I}$,
and an integer matrix $G'=(g'_{ij})_{i,j\in I}$
such that the following formulas hold:
\begin{align}
\label{eq:yfact}
y'_i &= \left(\prod_{j\in I} y_j^{c'_{ji}}\right)
\prod_{j\in I} F'_j(y)^{b'_{ji}}_{\oplus},
\\ 
\label{eq:xfact}
x'_i &= \left(\prod_{j\in I} x_j^{g'_{ji}}\right)
\frac{F'_i(\hat{y})}
{F'_i(y)_{\oplus}}.
\end{align}
Here, $F'_i(y)_{\oplus}$
is the sum  replacing the sum $+$
in any subtraction-free expression
of the polynomial $F'_i(y)$ with the sum
$\oplus$ in $\mathbb{P}_{\mathrm{univ}}(y)$.
\end{theorem}

We call the above $F'_i(y)$, $C'$, and $G'$ the {\em $F$-polynomials},
the {\em $C$-matrix}, and {\em $G$-matrix\/} of $(B',x',y')$, respectively.
Columns of $C'$ and  $G'$
are called {\em $\mathbf{c}$-vectors} and  {\em $\mathbf{g}$-vectors},
respectively, in \cite{Fomin07}.
It is conjectured in \cite[Sec.3]{Fomin02} that
the polynomial $F'_i(y)$ itself is subtraction-free (positivity conjecture).
However, we do not rely on this conjecture in this note.

The above formulas, together with categorification and tropicalization
explained later,
make the theory of cluster algebras very rich and  powerful.

\subsection{Sign-coherence and constant term conjectures}
\label{subsec:sign}

Let us recall two fundamental conjectures on cluster algebras,
described as `tantalizing conjecture' in \cite{Fomin07}.

\begin{conjecture}[{\cite[Conj.5.4]{Fomin07}}]
\label{conj:sign}
(i) (Sign-coherence conjecture)
Every column of a $C$-matrix
(i.e., every $\mathbf{c}$-vector)
is a nonzero vector and its nonzero components are either all positive
or all negative.
\par
(ii) (Constant term conjecture) Every $F$-polynomial has constant term 1.
\end{conjecture}

It is known that the conjectures (i) and (ii) are equivalent to each other  \cite[Prop.5.6]{Fomin07}.
Note that the conjecture (i) is analogous to the positivity/negativity property
of roots in root systems.
So far, the following partial result is known.

\begin{theorem}[\cite{Derksen10,Plamondon10a,Plamondon10b,Nagao10}]
The conjecture is true for any cluster algebra
$\mathcal{A}(B,x,y)$ with skew-symmetric matrix $B$.
\end{theorem}

All the known proofs of this theorem  require some `extra' machinery,
namely, representation of quiver with potential \cite{Derksen10},
categorification by 2-Calabi-Yau category \cite{Plamondon10a,Plamondon10b}
or by 3-Calabi-Yau category \cite{Nagao10}.

\subsection{Tropical sign}
Assuming the sign-coherence in Conjecture \ref{conj:sign},
to any seed $(B',x',y')$ and $k\in I$,
one can uniquely assign the sign $\varepsilon'_k\in
\{1,-1\}$ defined by the sign
of the $k$th column of the $C$-matrix of $(B',x',y')$.
We call the sign $\varepsilon'_k$ the {\em tropical sign
of $(B',x',y')$ at $k\in I$}.

In the next section,
we see that the tropical sign $\varepsilon'_k$ is the canonical choice
for $\varepsilon$ in the $\varepsilon$-expressions
\eqref{eq:bmut2}--\eqref{eq:xmut2}.
The use of $\varepsilon$-expressions, together with the tropical sign,
is important but overlooked until recently
\cite{Nagao10,Qin10,Nakanishi11a,Keller11,Nagao11b,Nagao11,Kashaev11}.

\section{Tropicalization in cluster algebras}

Here, we collect some 
basic properties of the tropicalization in cluster algebras.

\subsection{Tropical $y$-variables}

Let $(B',x',y')$ be any seed of a cluster algebra
$\mathcal{A}(B,x,y)$.
Let $\pi_{\mathrm{trop}}$ be the one
in \eqref{eq:hom}.
For the notational simplicity, let us write
$\pi_{\mathrm{trop}}(\,\cdot\,)$ as $[\,\cdot\,]$.
For the $y$-variables $y'_i \in \mathbb{P}_{\mathrm{univ}}(y)$,
we call their images $[y'_i]\in \mathbb{P}_{\mathrm{trop}}(y)$
the {\em tropical $y$-variables}.
They are  identified with the {\em principal coefficients} in
\cite{Fomin07}.

\begin{theorem}[{\cite[Prop.3.13 \& Prop.5.2]{Fomin07}}]
\begin{align}
\label{eq:ytrop}
[y'_i]&=\prod_{j\in I} y_j^{c'_{ji}},\\
\label{eq:f1}
[F'_i(y)_{\oplus}]&=1.
\end{align}
\end{theorem}
In other words, the formula \eqref{eq:yfact}
expresses the separation of the {\em tropical part
$[y'_i]$} and
the {\em nontropical part\/} given by the $F$-polynomials.
The sign-coherence conjecture
in  Conjecture \ref{conj:sign}
 means
that $[y'_i]$ is not 1 and its nonzero exponents in $y$
are either {\em all positive\/} or
{\em all negative}.
The constant term conjecture further implies that
$y'_i$ has a Laurent expansion in $y$,
and the tropical $y$-variable $[y'_i]$ is its {\em leading term}.
Note that \eqref{eq:f1} is weaker than
the constant term conjecture.
For example, $[y_1\oplus y_2]=1$.

Let us give the {\em exchange relation\/} for tropical
$y$-variables.
There are actually two versions,
 depending on
whether
one assumes the sign-coherence or not.

{\em (i) Without assuming the sign-coherence.}
First, note that the exchange relation of tropical $y$-variables
is formally the same with \eqref{eq:ymut2} by replacing
$y'_i$ with $[y'_i]$ and $\oplus$ with the
tropical sum $\oplus$ in \eqref{eq:tropo}.
Thus, we obtain the well-known exchange relation for $C$-matrices from
\eqref{eq:ymut2} \cite[Prop.5.8]{Fomin02}:
\begin{align}
\label{eq:cmut1}
c''_{ij} = 
\begin{cases}
-c'_{ik} & j=k\\
c'_{ij}+c'_{ik}[\varepsilon b'_{kj}]_+
+[-\varepsilon c'_{ik}]_+b'_{kj}
& j\neq k.
\end{cases}
\end{align}

{\em (ii) Assuming the sign-coherence.}
Next, assume the sign-coherence,
and let us take the tropical sign $\varepsilon'_k$
in Section \ref{subsec:sign}
and set $\varepsilon=\varepsilon'_k$
in \eqref{eq:ymut2}.
Then, we have
\begin{align}
[1\oplus y'_k{}^{\varepsilon'_k}]=1,
\end{align}
and we obtain a simplified form of the exchange
relation for tropical $y$-variables
\begin{align}
\label{eq:ymut3}
[y''_i]=
\begin{cases}
[y'_k]{}^{-1} & i= k\\
\displaystyle
[y'_i][ y'_k]{}^{[\varepsilon'_k b'_{ki}]_+}
& i\neq k,\\
\end{cases}
\end{align}
In other words,
by setting $\varepsilon=\varepsilon'_k$,
the second line of the exchange relation \eqref{eq:ymut2}
separates into
the tropical part
$y'_i y'_k{}^{[\varepsilon'_k b'_{ki}]_+}$
and the nontropical part
$(1\oplus y'_k{}^{\varepsilon'_k})^{-b'_{ki}}$.

Note that  \eqref{eq:ymut3} is equivalent to the exchange relation
of $C$-matrices
\begin{align}
\label{eq:cmut4}
c''_{ij} = 
\begin{cases}
-c'_{ik} & j=k\\
c'_{ij}+c'_{ik}[\varepsilon'_k b'_{kj}]_+
& j\neq k,
\end{cases}
\end{align}
which is also obtained directly from \eqref{eq:cmut1}  by
setting $\varepsilon=\varepsilon'_k$
 \cite[Prop.1.3]{Nakanishi11a}.

\subsection{$y$-tropical $x$-variables}
\label{subsec:ytropx}
In the same spirit, let us interpret the $G$-matrices as $x$-variables
in tropicalization.
However, it is {\em not\/} the  `tropical $x$-variables' in the direct
sense, {\em but\/} the one obtained through the tropicalization
 of $y$-variables (and also $\hat{y}$-variables!).
 Hence, we call them the 
 {\em $y$-tropical $x$-variables} here.
Again, we consider two versions.

{\em (i) Without assuming the sign-coherence.}
Let us recall the result of \cite{Fomin07}.
We introduce the
$\mathbb{Z}^n$-grading for the initial variables as follows:
\begin{align}
\label{eq:deg1}
\deg x_i = \vec{e}_i,
\quad
\deg y_i = -\vec{b}_i,
\end{align}
where $\vec{e}_i$ is the $i$th unit vector in 
$\mathbb{Z}^n$,
and
$\vec{b}_i$ is the $i$th column of $B$.
It follows that
\begin{align}
\label{eq:deg2}
\deg \hat{y}_i = \vec{0}.
\end{align}
Now we tropicalize the coefficients $y'_i \mapsto
[y'_i]$.
This also induces the evaluations
 $x'_i \mapsto [x'_i]$ and $\hat{y}'_i \mapsto [\hat{y}'_i]$,
where the coefficient of each monomial in $x$ is tropicalized.
More explicitly, from \eqref{eq:xfact} and \eqref{eq:f1},
we have
\begin{align}
\label{eq:xfact2}
[x'_i] &= \left(\prod_{j\in I} x_j^{g'_{ji}}\right)
F'_i(\hat{y}).
\end{align}
In particular, by \eqref{eq:deg1} and \eqref{eq:deg2},
$[x'_i]$ is homogeneous and
\begin{align}
\label{eq:xfact3}
\deg [x'_i] &= \vec{g}'_i,
\end{align}
where $\vec{g}'_i$ is the $i$th column of $G'$.
On the other hand, by the tropicalization of \eqref{eq:xmut2},
we have
\begin{align}
\label{eq:xmut3}
[x''_i]=
\begin{cases}
\displaystyle
[x'_k]{}^{-1}
\left(
\prod_{j\in I} [x'_j]{}^{[-\varepsilon b'_{jk}]_+}
\right)
\frac{1+[\hat{y}'_k]^{\varepsilon}}{[1\oplus y'_k{}^{\varepsilon}]}
& i=k\\
[x'_i] & i\neq k.\\
\end{cases}
\end{align}
Then, comparing
the degrees of the both hand sides of \eqref{eq:xmut3}, we have the exchange
relation of $G$-matrices \cite[Prop 6.6 \& Eq.(6.13)]{Fomin07}:
\begin{align}
\label{eq:gmut1}
g''_{ij}=
\begin{cases}
\displaystyle
-g'_{ik}
+
\sum_{\ell \in I}
g'_{i\ell}[-\varepsilon b'_{\ell k}]_+
-
\sum_{\ell \in I}
b_{i\ell}[-\varepsilon c'_{\ell k}]_+
& j = k\\
g'_{ij}
& j\neq k.\\
\end{cases}
\end{align}

{\em (ii) Assuming the sign-coherence.}
Next, assume the sign-coherence,
and let us take the tropical sign $\varepsilon'_k$
in Section \ref{subsec:sign}
and set $\varepsilon=\varepsilon'_k$
in \eqref{eq:xmut3}.
Then, the factor
${[1\oplus y'_k{}^{\varepsilon'_k}]}$ disappears.
One can now think that $[x_i]=x_i$ and $[\hat{y}_i]=\hat{y}_i$ are independent variables
and that 
the exchange relations for $[x'_i]$ and $[\hat{y}'_i]$
are given by
\eqref{eq:yhatmut} (with $\hat{y}'_i$ replaced by
$[\hat{y}'_i]$)  and \eqref{eq:xmut3}
(with $[1\oplus y'_k{}^{\varepsilon}]$
omitted).

We then make the `second tropicalization'
$[[\, \cdot \,]]:\mathbb{P}_{\mathrm{univ}}(\hat{y})
\rightarrow \mathbb{P}_{\mathrm{trop}}(\hat{y})$,
which also induces the evaluation
$[x'_i] \mapsto [[x'_i]]$.
Let us call the images $[[x'_i]]$ the {\em $y$-tropical $x$-variables}.
Then, from \eqref{eq:xfact2} and \eqref{eq:xmut3}, we obtain
\begin{align}
\label{eq:xfact4}
[[x'_i]] &= \prod_{j\in I} x_j^{g'_{ji}},\\
\label{eq:xmut4}
[[x''_i]]&=
\begin{cases}
\displaystyle
[[x'_k]]{}^{-1}
\prod_{j\in I} [[x'_j]]{}^{[-\varepsilon'_k b'_{jk}]_+}
& i=k\\
[[x'_i]] & i\neq k,\\
\end{cases}
\end{align}
which are parallel to \eqref{eq:ytrop} and \eqref{eq:ymut3}.
Thus,
by setting $\varepsilon=\varepsilon'_k$,
the second line of the exchange relation \eqref{eq:ymut2}
separates into
the $y$-tropical part and the $y$-nontropical part.

Note that  \eqref{eq:xmut4} is equivalent to the exchange relation
of $G$-matrices
\begin{align}
\label{eq:gmut4}
g''_{ij}=
\begin{cases}
\displaystyle
-g'_{ik}
+
\sum_{\ell \in I}
g'_{i\ell}[-\varepsilon'_k b'_{\ell k}]_+
& j = k\\
g'_{ij}
& j\neq k,\\
\end{cases}
\end{align}
which is also obtained directly from \eqref{eq:gmut1} by
setting $\varepsilon=\varepsilon'_k$
\cite[Prop.1.3]{Nakanishi11a}.

\subsection{Matrix form of mutations and duality}

{\em In this subsection we assume the sign-coherence.}

Following \cite{Berenstein05},
let $J_k$ be the $n \times n$ 
diagonal matrix  whose diagonal entries
are all 1 except for the $k$th one which is $-1$.
For any $n\times n$ matrix, say, $A$,
let $[A]^{k\bullet}_+$ denote the
matrix whose entries are zero except for the
$k$th row
and the $(k,i)$th entry
is $[a_{ki}]_+$.
Similarly,
let $[A]^{\bullet k}_+$  denote the
matrix whose entries are zero except for the
$k$th column
and the $(i,k)$th entry
is $[a_{ik}]_+$.
For an $n\times n$  matrix $A$, $k\in I$, and $\varepsilon\in \{1,-1\}$,
define
\begin{align}
\label{eq:mdef}
P_{A,k,\varepsilon}=
J_k + [\varepsilon A]^{k\bullet }_+,
\quad
Q_{A,k,\varepsilon}=
J_k + [\varepsilon A]^{\bullet k}_+.
\end{align}

Then, the exchange relation \eqref{eq:bmut2} with $\varepsilon=
\varepsilon'_k$,
\eqref{eq:cmut4}, and \eqref{eq:gmut4} are written in the following
matrix form
(\cite[Eq.(3.1)]{Berenstein05},
\cite[Prop.1.3]{Nakanishi11a}):
\begin{align}
\label{eq:mat1a}
B''&= Q_{B',k,-\varepsilon'_k} B' P_{B',k,\varepsilon'_k},
\\
\label{eq:mat1b}
C''&= C'P_{B',k,\varepsilon'_k},
\quad
G''= G' Q_{B',k,-\varepsilon'_k}.
\end{align}

For $(B'',x'',y'')=\mu_k(B',x',y')$,
we have $\varepsilon''_k = -\varepsilon'_k$,
$b''_{ik}=-b'_{ik}$,
and $b''_{ki}=-b'_{ki}$,
thus, we have
\begin{align}
\label{eq:mat2a}
P_{B'',k,\varepsilon''_k}
&=P_{B',k,\varepsilon'_k},\quad
Q_{B'',k,-\varepsilon''_k}
=Q_{B',k,-\varepsilon'_k},
\\
\label{eq:mat2b}
(P_{B',k,\varepsilon'_k})^2
&=
(Q_{B',k,-\varepsilon'_k})^2 = I,
\end{align}
which agrees with the involution property of the mutation
 $\mu_k^2=\mathrm{id}$.
 
Also, using $-B'{}^{T}=DB'D^{-1}$, we obtain
\begin{align}
Q_{B',k,-\varepsilon'_k}=
D^{-1} (P_{B',k,\varepsilon'_k})^{T} D.
\end{align}
Thus,
\eqref{eq:mat1a} and \eqref{eq:mat1b} are rephrased
in a more uniform way as
\begin{align}
\label{eq:mat3a}
DB''&= (P_{B',k,\varepsilon'_k})^T DB' P_{B',k,\varepsilon'_k},\\
\label{eq:mat3b}
DC''&= DC'P_{B',k,\varepsilon'_k},
\quad
G''D^{-1}= G'D^{-1} (P_{B',k,\varepsilon'_k})^{T}.
\end{align}

The following duality of $C$- and $G$-matrices 
immediately follows from \eqref{eq:mat2b}, \eqref{eq:mat3b},
and the fact $C=G=I$ for the initial seed $(B,x,y)$.
\begin{proposition}\cite[Eq.(3.11)]{Nakanishi11a}
\label{prop:duality}
 Assuming the sign-coherence, we have
\begin{align}
(G'D^{-1})^T (DC') = I.
\end{align}
\end{proposition}
In other words, the column vectors of $DC'$ and $G'D^{-1}$ are
dual basis of $\mathbb{Q}^n$ to each other.
(When $B$ is skew-symmetric, they are dual basis of
$\mathbb{Z}^n$ to each other.)

\subsection{Basis changes on lattice}

Let   ${L}\simeq \mathbb{Z}^n$ be a lattice of rank $n$.
Following Fock and Goncharov \cite{Fock03,Fock07,Fock07b},
we assign 
a basis $\vec{v}'_1,\dots,\vec{v}'_n$ of $L$
to each seed $\Sigma'=(B',x',y')$ 
such that under the mutation
$\Sigma''=\mu_k(\Sigma')$
two bases are related by
\begin{align}
\label{eq:ttrans}
 \vec{v}''_i&=
\begin{cases}
-\vec{v}'_k & i = k\\
\vec{v}'_i + [\varepsilon'_k b'_{ki}]_+ \vec{v}'_k 
& i\neq k.
\end{cases}
\end{align} 
This induces the coordinate transformation
 $\tau_{\Sigma',\Sigma''}:\mathbb{Z}^n \rightarrow \mathbb{Z}^n$
whose matrix representation is 
given by
$P_{B',k,\varepsilon'_k}$ in \eqref{eq:mdef}.
To each seed $\Sigma'$
we assign a skew-symmetric bilinear
form $\langle \, \cdot, \, \cdot \,
\rangle_{\Sigma'}$ on $L$ by
 $\langle \vec{v}'_i,\vec{v}'_j\rangle_{\Sigma'}
:= d_ib'_{ij}$.
The following is an immediate consequence of
\eqref{eq:mat3a}.
\begin{proposition}{\cite[Lemma 1.7]{Fock03}}
For any seeds $\Sigma'$ and $\Sigma''$, we have 
$\langle \vec{a}, \vec{b}\rangle_{\Sigma'}
=
\langle \vec{a},\vec{b} \rangle_{\Sigma''}
$
$(\vec{a},\vec{b}\in L)$.
Thus, the form $\langle \, \cdot\,, \,\cdot\, \rangle_{\Sigma'}$
does not depend on $\Sigma'$.
\end{proposition}

Fock and Goncharov formulated cluster algebras
(and quantum cluster algebras)
starting from the basis change
\eqref{eq:ttrans}
\cite{Fock03,Fock07,Fock07b}
without the tropical sign $\varepsilon'_k $.
Modifying the formulation to include $\varepsilon'_k$
is straightforward, though we need to establish the sign-coherence
in advance.
See \cite[Eq.(2.3)]{Nagao10} for an interpretation
of \eqref{eq:ttrans}
in view of Ginzburg's differential graded algebra.

\subsection{Categorification by generalized cluster category}

{\em In this subsection we assume that $B$ is skew-symmetric.}

The categorification of cluster algebras
by triangulated categories (2-Calabi-Yau realization)
was started by \cite{Caldero06,Buan06} with {\em cluster categories}
when the initial quiver $Q$ is type $ADE$.
Later it has been gradually extended to general quivers
\cite{Dehy08,Fu10,Palu08,Amiot09,Keller08,Keller09,
Plamondon10a,Plamondon10b}.
Here we summarize the most general result of \cite{
Plamondon10a,Plamondon10b} without explaining detail.

For the quiver $Q$ corresponding to $B$,
using the principal extension $\tilde{Q}$ of $Q$
and a potential $W$ on it,
a triangulated category
 $\mathcal{C}=\mathcal{C}_{(\tilde{Q},W)}$,
 called {\em the generalized cluster category of $Q$},
 is defined.
Then, to each seed $(B',x',y')$ of the cluster algebra
$\mathcal{A}(B,x,y)$,
one can canonically assign a rigid object $T'=\bigoplus_{i\in I} T'_i$
in $\mathcal{C}$.

\begin{theorem}[\cite{Plamondon10a,Plamondon10b}]
\label{thm:cat}
Let $T=\bigoplus_{i\in I} T_i$ be the rigid object
assigned to the initial seed $(B,x,y)$.
Then, for any seed $(B',x',y')$, the following formulas hold.
\begin{align}
\tilde{Q}' & = \mbox{quiver of\/ $\mathrm{End}_{\mathcal{C}}(T')$},\\
c'_{ij}& = - \mathrm{ind}_{T'}(T_i[1])_j
=\mathrm{ind}^{\mathrm{op}}_{T'}(T_i)_j,\\
\label{eq:gcat}
g'_{ij}& =  \mathrm{ind}_{T}(T'_j)_i,\\
F'_i(y)&=\sum_{e\in \mathbb{Z}^{\tilde{I}}_{\geq 0}}
\chi(\mathrm{Gr}_e(\mathrm{Hom}_{\mathcal{C}}(T,T'_i[1])))
\prod_{j\in I}y_j^{e_j}.
\end{align}
Here, $\mathrm{Gr}_e(\, \cdot\, )$
the quiver Grassmannian with dimension vector $e$,
and $\chi(\, \cdot \, )$ is the Euler characteristic.
\end{theorem}

In brief, the rigid object $T'$ provides all the information
of the seed $(B',x',y')$.
Let us present some consequence of this remarkable theorem
to tropicalization.

Recall that the matrices $C'$ and $G'$ determine each other
(Proposition \ref{prop:duality}).
By \eqref{eq:gcat}, the matrix $G'$ carries the information
of the index of $T'$.
Moreover, the index of $ T'$ uniquely determines
$T'$ itself \cite[Prop.3.1]{Plamondon10b}.
Thus, by Theorem \ref{thm:cat}, we recover the seed $(B',x',y')$.
Therefore, we obtain the following corollary.

\begin{corollary}
The tropical $y$-variables $[y'_i]$ $(i\in I)$ determines
the seed $(B',x',y')$.
\end{corollary}

In our application, it is useful to formulate this
corollary in terms of periodicity.

\begin{definition}
\label{def:per}
Let $\nu: I \rightarrow I$ be a bijection,
and let $(B',x',y')$ be a seed.
For an $I$-sequence $(k_1,\dots,k_L)$, let
$(B'',x'',y'')=
\mu_{k_L} \cdots \mu_{k_1}(B',x',y')$.
We say $(k_1,\dots,k_L)$ is a {\em $\nu$-period of $(B',x',y')$}
if
\begin{align}
\label{eq:per}
&b''_{\nu(i)\nu(j)}=b'_{ij},\quad
x''_{\nu(i)}=x'_{i},\quad
y''_{\nu(i)}=y'_{i},\quad
(i,j\in I).
\end{align}
\end{definition}

\begin{example} In Example \ref{ex:a2}
$(1,2,1,2,1)$ is a $(12)$-period of $(Q,x,y)$,
where $(12)$ is the transposition of $1$ and $2$.
\end{example}

\begin{corollary}[{\cite[Th.5.1]{Inoue10a}}]
\label{cor:yperiod}
An $I$-sequence $(k_1,\dots,k_L)$ is a $\nu$-period of $(B',x',y')$
if and only if
$$
[y''_{\nu(i)}] =[y'_i] \quad (i\in I).
$$
\end{corollary}
In other words,
the periodicity of seeds follows from
the periodicity of tropical $y$-variables.
(We conjecture that it is true also for skew-symmetrizable $B$.)
It turns out that this criterion is very powerful.

\section{Application I: Periodicities of Y-systems}

We apply the tropicalization method to
prove certain periodicities of seeds
corresponding to the conjectured periodicities of Y-systems.

\subsection{Background: Periodicity conjecture of Y-systems}
\label{subsec:period}

For a cluster algebra $\mathcal{A}(B,x,y)$ of finite type,
starting from $(B,x,y)$, repeat any finite sequence of mutations 
sufficiently many times.
Since there are only finitely many numbers of seeds by definition,
eventually it returns to the initial quiver.
Thus, we have a periodicity of seeds.
Aside from these examples, however,
periodicity of seeds is a rare event
and so far there is no systematic way to find it.
Fortunately, we had some source  predicting
infinitely many periodicities of seeds,
 which appeared prior to cluster algebras;
that is, the {\em periodicity conjecture of Y-systems}
\cite{Zamolodchikov91,Ravanini93,Kuniba94a}.

The Y-systems are systems of functional algebraic equations,
and they were introduced and studied in 90's
(see \cite{Kuniba10} and references therein).
They are associated with pairs $(X,\ell)$ of
a Dynkin diagram $X$ of finite type
 and an integer $\ell\geq 2$ (called {\em level\/}).
For simplicity, let us concentrate on the case when $X$ is simply laced. Let $\mathcal{I}$ be the set of vertices of $X$,
and consider a family of variables
$\{ Y^{a}_m(u) \mid
a\in \mathcal{I};
m=1,\dots,\ell-1; u\in \mathbb{Z}\}$,
The Y-system is the following system of equations:
\begin{align}
\label{eq:Ysystem}
Y^{(a)}_m(u-1)Y^{(a)}_m(u+1)
=
\frac{
\prod_{b\in \mathcal{I}: b\sim a} (1+Y^{(b)}_m(u))
}
{
(1+Y^{(a)}_{m-1}(u)^{-1})(1+Y^{(a)}_{m+1}(u)^{-1})
}.
\end{align}
Here,  $a\sim b$ means
$a$ and $b$ are adjacent in $X$,
and $Y^{(a)}_{0}(u)^{-1}=Y^{(a)}_{\ell}(u)^{-1}=0$.
The periodicity conjecture of the Y-system claimed that
\begin{align}
\label{eq:yperiod}
Y^{(a)}_m(u+2(h+\ell))
= Y^{(a)}_m(u),
\end{align}
where $h$ is the Coxeter number of $X$.
When $X$ is nonsimply laced, the Y-system looks more complicated
and we omit to present it,
but the periodicity \eqref{eq:yperiod} still holds
by replacing $h$ with the dual Coxeter number $h^{\vee}$.
(For simply laced case $h=h^{\vee}$.)
The periodicity \eqref{eq:yperiod} was proved
for $X=A_1$ by \cite{Frenkel95, Gliozzi96} 
and 
for $X=A_n$ by \cite{Szenes09,Volkov07},
by using explicit solutions of 
\eqref{eq:Ysystem}.

It turned out that these $Y$-systems are a part of
relations among $y$-variables in certain cluster algebras,
and the periodicity conjecture
\eqref{eq:yperiod}
is translated into the periodicity of seeds
in the corresponding cluster algebras.
Fomin and Zelevinsky first
recognized this fact in the simplest case
of simply laced $X$ and $\ell=2$,
and proved the periodicity  \cite{Fomin03b}.
Then, the case  of
simply laced $X$ and {\em general} $\ell$
was proved using the categorification
and the Auslander-Reiten theory of quiver representations
\cite{Keller08,Keller10}.
Finally, the most general case including the nonsimply laced $X$
was proved using the tropicalization method.
\cite{Inoue10a,Inoue10b}.
Below we explain the essence of this method.


\subsection{Sign-arrow coordination}
\label{subsec:sa}
Fortunately, for the cluster algebra corresponding
to any Y-system,
the initial exchange matrix $B$ is   skew-symmetric,
even when $X$ is nonsimply laced.
Thus, one can apply Corollary \ref{cor:yperiod}
and the sign-coherence.

To study periodicity, it is efficient to work with
the quiver $Q$  corresponding to $B$.
Recall the exchange relation  \eqref{eq:ymut3} of tropical
$y$-variables.
The only nontrivial term therein is
$[ y'_k]{}^{[\varepsilon'_k b'_{ki}]_+}$
and this term appears only when $\varepsilon'_k b'_{ki}>0$.
In particular,
when there is no multiple arrows between
$k$ and $i$,
which is the case for the Y-systems,
the tropical exchange relation \eqref{eq:ymut3}
is rephrased in the following simple rule:

\begin{itemize}
\item[i)] (Sign-arrow coordination) For $i\neq k$,
multiply $y'_k$ to $y'_i$ only when one of the
following situations occurs:
\begin{align}
\label{eq:signarrow}
\raisebox{0pt}{
\begin{picture}(200,5)(0,-5)
\put(0,0){{{\circle{8}}}}
\put(0,0){{\circle{3}}}
\put(30,0){{\circle{3}}}
\put(3,0){\vector(1,0){24}}
\put(-3,-12){\small $k$}
\put(27,-12){\small $i$}
\put(50,-3){$\varepsilon'_k>0$}
\put(100,-3){or}
\put(130,0)
{
\put(0,0){{{\circle{8}}}}
\put(0,0){{\circle{3}}}
\put(30,0){{\circle{3}}}
\put(27,0){\vector(-1,0){24}}
\put(-3,-12){\small $k$}
\put(27,-12){\small $i$}
\put(50,-3){$\varepsilon'_k<0$}
}
\end{picture}
}
\end{align}

\item[ii)]
Invert $y'_k$.
\end{itemize}

\subsection{Examples} Let us demonstrate,
by examples, how the tropicalization method
works to prove periodicities of seeds corresponding to
 \eqref{eq:yperiod}.
Our goal here is to understand the period $2(h+\ell)$ in  \eqref{eq:yperiod}
for the case $(X,\ell)=(A_3,3)$ in Example
\ref{ex:a3l3}.
Recall that the Coxeter number of $A_n$ is $n+1$.

\begin{example} {\em The case $(X,\ell)=(A_2,2)$,
where $h+\ell= 3+2 = 5$.}
\label{ex:a2l2}

The corresponding cluster algebra 
is of type $A_2$,
and the (half) periodicity of seeds is just the one presented in
Example \ref{ex:a2}.
The length 5 of the (half) period $(1,2,1,2,1)$ 
coincides with $5 = h + \ell=3+2$.

Using Corollary \ref{cor:yperiod}
and the rules i) and ii) in Section \ref{subsec:sa},
one can prove this periodicity {\em much more\/} efficiently
by the following diagrammatic calculation of tropical $y$-variables:
\begin{align*}
\begin{picture}(260,73)(-20,-10)
\put(0,40){\vector(0,-1){20}}
\put(-5, -3){\framebox(10,20)[c]}
\put(-2,9){$\scriptstyle 0$}
\put(-2,1){$\scriptstyle 1$}
\put(-2,55){$\scriptstyle 1$}
\put(-2,47){$\scriptstyle 0$}
\put(-8,-15){$y(0)$}
\put(-20,5){$y_1$}
\put(-20,50){$y_2$}
\put(15,25){$\leftrightarrow$}
\put(15,35){$\mu_1$}
\put(40,20){\vector(0,1){20}}
\put(38,9){$\scriptstyle 0$}
\put(36,1){-$\scriptstyle 1$}
\put(35, 43){\framebox(10,20)[c]}
\put(38,55){$\scriptstyle 1$}
\put(38,47){$\scriptstyle 0$}
\put(32,-15){$y(1)$}
\put(40,0)
{
\put(15,25){$\leftrightarrow$}
\put(15,35){$\mu_2$}
\put(40,40){\vector(0,-1){20}}
\put(35, -3){\framebox(10,20)[c]}
\put(38,9){$\scriptstyle 0$}
\put(36,1){-$\scriptstyle 1$}
\put(36,55){-$\scriptstyle 1$}
\put(38,47){$\scriptstyle 0$}
\put(32,-15){$y(2)$}
}
\put(80,0)
{
\put(15,25){$\leftrightarrow$}
\put(15,35){$\mu_1$}
\put(40,20){\vector(0,1){20}}
\put(38,9){$\scriptstyle 0$}
\put(38,1){$\scriptstyle 1$}
\put(35, 43){\framebox(10,20)[c]}
\put(36,55){-$\scriptstyle 1$}
\put(36,47){-$\scriptstyle 1$}
\put(32,-15){$y(3)$}
}
\put(120,0)
{
\put(15,25){$\leftrightarrow$}
\put(15,35){$\mu_2$}
\put(40,40){\vector(0,-1){20}}
\put(35, -3){\framebox(10,20)[c]}
\put(36,9){-$\scriptstyle 1$}
\put(38,1){$\scriptstyle 0$}
\put(38,55){$\scriptstyle 1$}
\put(38,47){$\scriptstyle 1$}
\put(32,-15){$y(4)$}
}
\put(160,0)
{
\put(15,25){$\leftrightarrow$}
\put(15,35){$\mu_1$}
\put(40,20){\vector(0,1){20}}
\put(35, 43){\framebox(10,20)[c]}
\put(38,9){$\scriptstyle 1$}
\put(38,1){$\scriptstyle 0$}
\put(38,55){$\scriptstyle 0$}
\put(38,47){$\scriptstyle 1$}
\put(32,-15){$y(5)$}
}
\end{picture}
\end{align*}
Here,
the vectors $\left(\begin{matrix}0\\1\end{matrix}\right)$
and $\left(\begin{matrix}-1\\-1\end{matrix}\right)$,
for example, represent 
$y_1$ and $y_1^{-1} y_2^{-1}$, respectively.
Thus, they are $\mathbf{c}$-vectors (leading from the bottom).
Observe that, at the forward mutation points (which are framed in the diagram),
the exponents of the variables $[y_i(u)]$ for
$u=0,1$ are {\em positive} and identified with
the {\em negative simple roots\/} $-\alpha_1$, $-\alpha_2$ of $A_2$,
while the ones for $u=2,3,4$
are {\em negative} and identified with
the {\em positive roots\/} $\alpha_1$, $\alpha_1+\alpha_2$,
$\alpha_2$ of $A_2$.
This separation of the region of $u$ naturally
explains the formula of the period  $5=3+2$.
Fomin and Zelevinsky explained and generalized this
phenomenon to any simply laced $X$,
using the piecewise-linear analogue of the Coxeter element  \cite{Fomin03b}.

For later use, let us present a similar diagram for the opposite quiver of $Q$.
\begin{align*}
\begin{picture}(255,73)(-15,-10)
\put(0,20){\vector(0,1){20}}
\put(-5, -3){\framebox(10,20)[c]}
\put(-2,9){$\scriptstyle 0$}
\put(-2,1){$\scriptstyle 1$}
\put(-2,55){$\scriptstyle 1$}
\put(-2,47){$\scriptstyle 0$}
\put(-8,-15){$y(0)$}
\put(15,25){$\leftrightarrow$}
\put(15,35){$\mu_1$}
\put(40,40){\vector(0,-1){20}}
\put(38,9){$\scriptstyle 0$}
\put(36,1){-$\scriptstyle 1$}
\put(35, 43){\framebox(10,20)[c]}
\put(38,55){$\scriptstyle 1$}
\put(38,47){$\scriptstyle 1$}
\put(32,-15){$y(1)$}
\put(40,0)
{
\put(15,25){$\leftrightarrow$}
\put(15,35){$\mu_2$}
\put(40,20){\vector(0,1){20}}
\put(35, -3){\framebox(10,20)[c]}
\put(38,9){$\scriptstyle 1$}
\put(38,1){$\scriptstyle 0$}
\put(36,55){-$\scriptstyle 1$}
\put(36,47){-$\scriptstyle 1$}
\put(32,-15){$y(2)$}
}
\put(80,0)
{
\put(15,25){$\leftrightarrow$}
\put(15,35){$\mu_1$}
\put(40,40){\vector(0,-1){20}}
\put(36,9){-$\scriptstyle 1$}
\put(38,1){$\scriptstyle 0$}
\put(35, 43){\framebox(10,20)[c]}
\put(38,55){$\scriptstyle  0$}
\put(36,47){-$\scriptstyle 1$}
\put(32,-15){$y(3)$}
}
\put(120,0)
{
\put(15,25){$\leftrightarrow$}
\put(15,35){$\mu_2$}
\put(40,20){\vector(0,1){20}}
\put(35, -3){\framebox(10,20)[c]}
\put(36,9){-$\scriptstyle 1$}
\put(38,1){$\scriptstyle 0$}
\put(38,55){$\scriptstyle 0$}
\put(38,47){$\scriptstyle 1$}
\put(32,-15){$y(4)$}
}
\put(160,0)
{
\put(15,25){$\leftrightarrow$}
\put(15,35){$\mu_1$}
\put(40,40){\vector(0,-1){20}}
\put(35, 43){\framebox(10,20)[c]}
\put(38,9){$\scriptstyle 1$}
\put(38,1){$\scriptstyle 0$}
\put(38,55){$\scriptstyle 0$}
\put(38,47){$\scriptstyle 1$}
\put(32,-15){$y(5)$}
}
\end{picture}
\end{align*}
\end{example}

\begin{example} {\em The case $(X,\ell)=(A_3,2)$,
where $h+\ell= 4+2 = 6$.}
\label{ex:a3l2}

Let us present one more example of
level $2$.
The corresponding cluster algebra
is the cluster algebra of type $A_3$.
We choose the following initial quiver $Q$
with parity $\pm$ assignment to each vertex.
\begin{align}
Q =\ 
\raisebox{-10pt}{
\begin{picture}(80,25)(0,-15)
%
%
\put(0,0){\circle{6}}
\put(40,0){\circle{6}}
\put(80,0){\circle{6}}
\put(37,0){\vector(-1,0){34}}
\put(43,0){\vector(1,0){34}}
\put(-3,-15){$1$}
\put(37,-15){$2$}
\put(77,-15){$3$}
\put(-3,8){$+$}
\put(37,8){$-$}
\put(77,8){$+$}
\end{picture}
}
\, .
\end{align}
Let $\mu_+=\mu_1\mu_3$ and $\mu_-=\mu_2$.
Set $\Sigma(0)=(Q(0),x(0),y(0))=(Q,x,y)$,
and consider the following sequence of mutations
in the {\em backward direction}:
\begin{align}
\label{eq:seedmutseq5}
\Sigma(-6)
\
\mathop{\leftrightarrow}^{\mu_{+}}
\
\Sigma(-5)
\
\mathop{\leftrightarrow}^{\mu_{-}}
\
\Sigma(-4)
\
\mathop{\leftrightarrow}^{\mu_{+}}
\
\Sigma(-3)
\
\mathop{\leftrightarrow}^{\mu_{-}}
\
\Sigma(-2)
\
\mathop{\leftrightarrow}^{\mu_{+}}
\
\Sigma(-1)
\
\mathop{\leftrightarrow}^{\mu_{-}}
\
\Sigma(0)
\end{align}
By doing the backward mutations
(at the points  which are {\em not\/} framed in the diagram below)
for tropical $y$-variables,
one can show the (half) periodicity.
\begin{align*}
\begin{picture}(253,110)(-15,18)
\put(0,120)
{
\put(-3, -5){\framebox(26,10)[c]}
\put(0,-2){$\scriptstyle 0$}
\put(8,-2){$\scriptstyle 0$}
\put(16,-2){$\scriptstyle  1$}
\put(35,0){\vector(-1,0){10}}
\put(40,-2){$\scriptstyle 0$}
\put(48,-2){$\scriptstyle 1$}
\put(56,-2){$\scriptstyle 0$}
\put(65,0){\vector(1,0){10}}
\put(77, -5){\framebox(26,10)[c]}
\put(80,-2){$\scriptstyle 1$}
\put(88,-2){$\scriptstyle 0$}
\put(96,-2){$\scriptstyle 0$}
\put(40,-15){$y(-6)$}
}
\put(150,120)
{
\put(37, -5){\framebox(26,10)[c]}
\put(0,-2){$\scriptstyle 0$}
\put(8,-2){$\scriptstyle 0$}
\put(14,-2){-$\scriptstyle  1$}
\put(25,0){\vector(1,0){10}}
\put(40,-2){$\scriptstyle 0$}
\put(48,-2){$\scriptstyle 1$}
\put(56,-2){$\scriptstyle 0$}
\put(75,0){\vector(-1,0){10}}
\put(78,-2){-$\scriptstyle 1$}
\put(88,-2){$\scriptstyle 0$}
\put(96,-2){$\scriptstyle 0$}
\put(40,-15){$y(-5)$}
}
\put(120,115){$\leftrightarrow$}
\put(120,125){$\mu_+$}
\put(0,90)
{
\put(-3, -5){\framebox(26,10)[c]}
\put(0,-2){$\scriptstyle 0$}
\put(8,-2){$\scriptstyle 0$}
\put(14,-2){-$\scriptstyle  1$}
\put(35,0){\vector(-1,0){10}}
\put(40,-2){$\scriptstyle 0$}
\put(46,-2){-$\scriptstyle 1$}
\put(56,-2){$\scriptstyle 0$}
\put(65,0){\vector(1,0){10}}
\put(77, -5){\framebox(26,10)[c]}
\put(78,-2){-$\scriptstyle 1$}
\put(88,-2){$\scriptstyle 0$}
\put(96,-2){$\scriptstyle 0$}
\put(40,-15){$y(-4)$}
}
\put(-30,85){$\leftrightarrow$}
\put(-30,95){$\mu_-$}
\put(150,90)
{
\put(37, -5){\framebox(26,10)[c]}
\put(0,-2){$\scriptstyle 0$}
\put(8,-2){$\scriptstyle 0$}
\put(16,-2){$\scriptstyle  1$}
\put(25,0){\vector(1,0){10}}
\put(38,-2){-$\scriptstyle 1$}
\put(46,-2){-$\scriptstyle 1$}
\put(54,-2){-$\scriptstyle 1$}
\put(75,0){\vector(-1,0){10}}
\put(80,-2){$\scriptstyle 1$}
\put(88,-2){$\scriptstyle 0$}
\put(96,-2){$\scriptstyle 0$}
\put(40,-15){$y(-3)$}
}
\put(120,85){$\leftrightarrow$}
\put(120,95){$\mu_+$}
\put(0,60)
{
\put(-3, -5){\framebox(26,10)[c]}
\put(-2,-2){-$\scriptstyle 1$}
\put(6,-2){-$\scriptstyle 1$}
\put(16,-2){$\scriptstyle  0$}
\put(35,0){\vector(-1,0){10}}
\put(40,-2){$\scriptstyle 1$}
\put(48,-2){$\scriptstyle 1$}
\put(56,-2){$\scriptstyle 1$}
\put(65,0){\vector(1,0){10}}
\put(77, -5){\framebox(26,10)[c]}
\put(80,-2){$\scriptstyle 0$}
\put(86,-2){-$\scriptstyle 1$}
\put(94,-2){-$\scriptstyle 1$}
\put(40,-15){$y(-2)$}
}
\put(-30,55){$\leftrightarrow$}
\put(-30,65){$\mu_-$}
\put(150,60)
{
\put(37, -5){\framebox(26,10)[c]}
\put(0,-2){$\scriptstyle 1$}
\put(8,-2){$\scriptstyle 1$}
\put(16,-2){$\scriptstyle  0$}
\put(25,0){\vector(1,0){10}}
\put(40,-2){$\scriptstyle 0$}
\put(46,-2){-$\scriptstyle 1$}
\put(56,-2){$\scriptstyle 0$}
\put(75,0){\vector(-1,0){10}}
\put(80,-2){$\scriptstyle 0$}
\put(88,-2){$\scriptstyle 1$}
\put(96,-2){$\scriptstyle 1$}
\put(40,-15){$y(-1)$}
}
\put(120,55){$\leftrightarrow$}
\put(120,65){$\mu_+$}
\put(0,30)
{
\put(-3, -5){\framebox(26,10)[c]}
\put(0,-2){$\scriptstyle 1$}
\put(8,-2){$\scriptstyle 0$}
\put(16,-2){$\scriptstyle  0$}
\put(35,0){\vector(-1,0){10}}
\put(40,-2){$\scriptstyle 0$}
\put(48,-2){$\scriptstyle 1$}
\put(56,-2){$\scriptstyle 0$}
\put(65,0){\vector(1,0){10}}
\put(77, -5){\framebox(26,10)[c]}
\put(80,-2){$\scriptstyle 0$}
\put(88,-2){$\scriptstyle 0$}
\put(96,-2){$\scriptstyle 1$}
\put(40,-15){$y(0)$}
}
\put(-30,25){$\leftrightarrow$}
\put(-30,35){$\mu_-$}
\end{picture}
\end{align*}
Here, the vectors $(1,0,0)$ and $(0,-1,-1)$, for example,
represent
$y_1$ and $y_2^{-2}y_3^{-3}$, respectively:
Again, at the forward mutation points,
(which are framed in the diagram), $[y_i(u)]$ for
$u=-1,-2,-3,-4$ are identified with
the {\em positive roots\/}
$\alpha_1$, $\alpha_2$,  $\alpha_3$, $\alpha_1+\alpha_2$,
$\alpha_2+\alpha_3$, $\alpha_1+\alpha_2+\alpha_3$,
of $A_3$,
while
the ones for $u=-5,-6$ are identified with
the {\em negative simple roots\/} $-\alpha_1$,
$-\alpha_2$, $-\alpha_3$ of $A_3$.
This explains the formula for the (half) period $6=4+2$.

For later use, let us present
a similar diagram for the opposite quiver
with opposite parity:
\begin{align*}
\begin{picture}(253,110)(-15,18)
\put(0,120)
{
\put(0,-2){$\scriptstyle 0$}
\put(8,-2){$\scriptstyle 0$}
\put(16,-2){$\scriptstyle  1$}
\put(25,0){\vector(1,0){10}}
\put(37, -5){\framebox(26,10)[c]}
\put(40,-2){$\scriptstyle 0$}
\put(48,-2){$\scriptstyle 1$}
\put(56,-2){$\scriptstyle 0$}
\put(75,0){\vector(-1,0){10}}
%
\put(80,-2){$\scriptstyle 1$}
\put(88,-2){$\scriptstyle 0$}
\put(96,-2){$\scriptstyle 0$}
\put(40,-15){$y(-6)$}
}
\put(150,120)
{
\put(-3, -5){\framebox(26,10)[c]}
\put(0,-2){$\scriptstyle 0$}
\put(8,-2){$\scriptstyle 0$}
\put(16,-2){$\scriptstyle  1$}
\put(35,0){\vector(-1,0){10}}
%
\put(40,-2){$\scriptstyle 0$}
\put(46,-2){-$\scriptstyle 1$}
\put(56,-2){$\scriptstyle 0$}
\put(65,0){\vector(1,0){10}}
\put(77, -5){\framebox(26,10)[c]}
\put(80,-2){$\scriptstyle 1$}
\put(88,-2){$\scriptstyle 0$}
\put(96,-2){$\scriptstyle 0$}
\put(40,-15){$y(-5)$}
}
\put(120,115){$\leftrightarrow$}
\put(120,125){$\mu_+$}
\put(0,90)
{
\put(0,-2){$\scriptstyle 0$}
\put(8,-2){$\scriptstyle 0$}
\put(14,-2){-$\scriptstyle  1$}
\put(25,0){\vector(1,0){10}}
\put(37, -5){\framebox(26,10)[c]}
\put(40,-2){$\scriptstyle 0$}
\put(46,-2){-$\scriptstyle 1$}
\put(56,-2){$\scriptstyle 0$}
\put(75,0){\vector(-1,0){10}}
%
\put(78,-2){-$\scriptstyle 1$}
\put(88,-2){$\scriptstyle 0$}
\put(96,-2){$\scriptstyle 0$}
\put(40,-15){$y(-4)$}
}
\put(-30,85){$\leftrightarrow$}
\put(-30,95){$\mu_-$}
\put(150,90)
{
\put(-3, -5){\framebox(26,10)[c]}
\put(0,-2){$\scriptstyle 0$}
\put(6,-2){-$\scriptstyle 1$}
\put(14,-2){-$\scriptstyle  1$}
\put(35,0){\vector(-1,0){10}}
%
\put(40,-2){$\scriptstyle 0$}
\put(48,-2){$\scriptstyle 1$}
\put(56,-2){$\scriptstyle 0$}
\put(65,0){\vector(1,0){10}}
\put(77, -5){\framebox(26,10)[c]}
\put(78,-2){-$\scriptstyle 1$}
\put(86,-2){-$\scriptstyle 1$}
\put(96,-2){$\scriptstyle 0$}
\put(40,-15){$y(-3)$}
}
\put(120,85){$\leftrightarrow$}
\put(120,95){$\mu_+$}
\put(0,60)
{
\put(0,-2){$\scriptstyle 0$}
\put(8,-2){$\scriptstyle 1$}
\put(16,-2){$\scriptstyle  1$}
\put(25,0){\vector(1,0){10}}
\put(37, -5){\framebox(26,10)[c]}
\put(38,-2){-$\scriptstyle 1$}
\put(46,-2){-$\scriptstyle 1$}
\put(54,-2){-$\scriptstyle 1$}
\put(75,0){\vector(-1,0){10}}
%
\put(80,-2){$\scriptstyle 1$}
\put(88,-2){$\scriptstyle 1$}
\put(96,-2){$\scriptstyle 0$}
\put(40,-15){$y(-2)$}
}
\put(-30,55){$\leftrightarrow$}
\put(-30,65){$\mu_-$}
\put(150,60)
{
\put(-2,-2){-$\scriptstyle 1$}
\put(8,-2){$\scriptstyle 0$}
\put(16,-2){$\scriptstyle  0$}
\put(35,0){\vector(-1,0){10}}
\put(-3, -5){\framebox(26,10)[c]}
\put(40,-2){$\scriptstyle 1$}
\put(48,-2){$\scriptstyle 1$}
\put(56,-2){$\scriptstyle 1$}
\put(65,0){\vector(1,0){10}}
\put(77, -5){\framebox(26,10)[c]}
\put(80,-2){$\scriptstyle 0$}
\put(88,-2){$\scriptstyle 0$}
\put(94,-2){-$\scriptstyle 1$}
\put(40,-15){$y(-1)$}
}
\put(120,55){$\leftrightarrow$}
\put(120,65){$\mu_+$}
\put(0,30)
{
\put(0,-2){$\scriptstyle 1$}
\put(8,-2){$\scriptstyle 0$}
\put(16,-2){$\scriptstyle  0$}
\put(25,0){\vector(1,0){10}}
\put(37, -5){\framebox(26,10)[c]}
\put(40,-2){$\scriptstyle 0$}
\put(48,-2){$\scriptstyle 1$}
\put(56,-2){$\scriptstyle 0$}
\put(75,0){\vector(-1,0){10}}
%
\put(80,-2){$\scriptstyle 0$}
\put(88,-2){$\scriptstyle 0$}
\put(96,-2){$\scriptstyle 1$}
\put(40,-15){$y(0)$}
}
\put(-30,25){$\leftrightarrow$}
\put(-30,35){$\mu_-$}
\end{picture}
\end{align*}

\end{example}

\begin{example}
 {\em The case $(X,\ell)=(A_3,3)$ with
$h+\ell= 4+3 = 7$.}
\label{ex:a3l3}

This is our main example.
The corresponding cluster algebra  has the following initial quiver $Q$
with parity  assignment to each vertex.
\begin{align}
Q =\ 
\quad
\raisebox{-27.5pt}{
\begin{picture}(80,55)(0,-10)
\put(0,0){\circle{6}}
\put(40,0){\circle{6}}
\put(80,0){\circle{6}}
\put(0,40){\circle{6}}
\put(40,40){\circle{6}}
\put(80,40){\circle{6}}
\put(37,0){\vector(-1,0){34}}
\put(43,0){\vector(1,0){34}}
\put(3,40){\vector(1,0){34}}
\put(77,40){\vector(-1,0){34}}
\put(0,3){\vector(0,1){34}}
\put(40,37){\vector(0,-1){34}}
\put(80,3){\vector(0,1){34}}
\put(-10,-10){$+$}
\put(30,-10){$-$}
\put(70,-10){$+$}
\put(-10,30){$-$}
\put(30,30){$+$}
\put(70,30){$-$}
\end{picture}
}
\, .
\end{align}
Let $\mu_+$ and $\mu_-$ be the composite mutations
at the vertices with $+$ and $-$, respectively.
Then, we consider the sequence of mutations
in the both forward and backward directions:
\begin{align}
\label{eq:seedmutseq6}
\cdots
\
\mathop{\leftrightarrow}^{\mu_{-}}
\
\Sigma(-2)
\
\mathop{\leftrightarrow}^{\mu_{+}}
\
\Sigma(-1)
\
\mathop{\leftrightarrow}^{\mu_{-}}
\
\Sigma(0)
\
\mathop{\leftrightarrow}^{\mu_{+}}
\
\Sigma(1)
\
\mathop{\leftrightarrow}^{\mu_{-}}
\
\Sigma(2)
\
\mathop{\leftrightarrow}^{\mu_{+}}
\
\cdots.
\end{align}
We can prove the desired (half) period of $7$
by the following diagrammatic calculation.
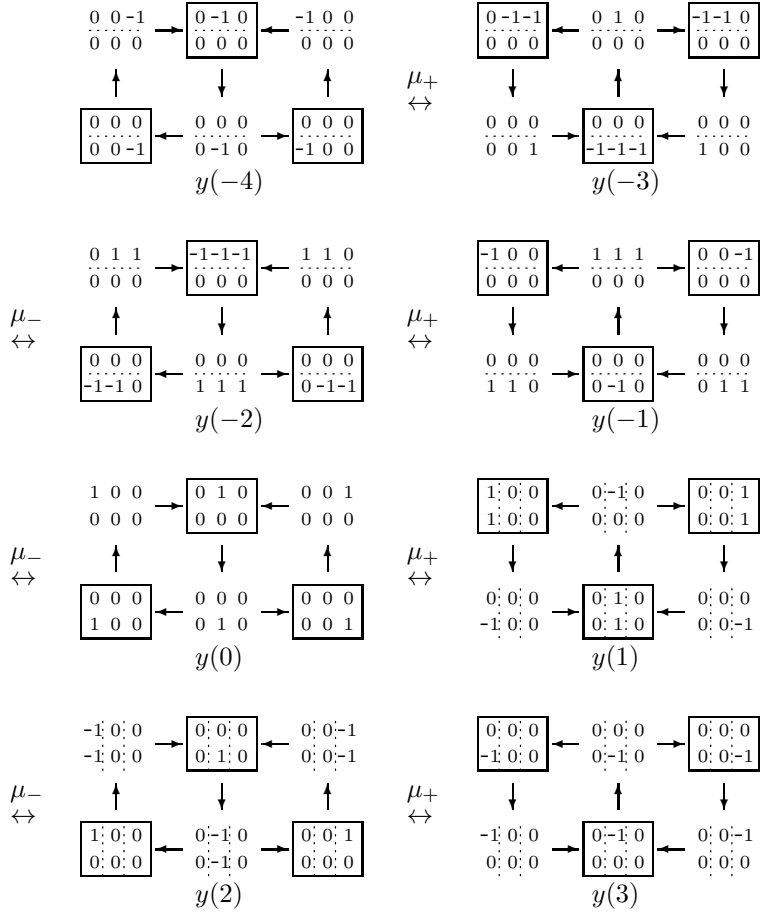
\begin{center}
\begin{picture}(255,345)(-15,-15)
\put(0,270)
{
%
\put(0,48){$\scriptstyle 0$}
\put(8,48){$\scriptstyle 0$}
\put(14,48){-$\scriptstyle  1$}
\put(0,38){$\scriptstyle 0$}
\put(8,38){$\scriptstyle 0$}
\put(16,38){$\scriptstyle  0$}
\dottedline{3}(0,45)(20,45)%
\put(10,20){\vector(0,1){10}}
\put(25,45){\vector(1,0){10}}
\put(37, 35){\framebox(26,20)[c]}
\put(40,48){$\scriptstyle 0$}
\put(46,48){-$\scriptstyle 1$}
\put(56,48){$\scriptstyle 0$}
\put(40,38){$\scriptstyle 0$}
\put(48,38){$\scriptstyle 0$}
\put(56,38){$\scriptstyle 0$}
\dottedline{3}(40,45)(60,45)
\put(50,30){\vector(0,-1){10}}
\put(75,45){\vector(-1,0){10}}
%
\put(78,48){-$\scriptstyle 1$}
\put(88,48){$\scriptstyle 0$}
\put(96,48){$\scriptstyle 0$}
\put(80,38){$\scriptstyle 0$}
\put(88,38){$\scriptstyle 0$}
\put(96,38){$\scriptstyle 0$}
\dottedline{3}(80,45)(100,45)
\put(90,20){\vector(0,1){10}}
\put(-3, -5){\framebox(26,20)[c]}
\put(0,8){$\scriptstyle 0$}
\put(8,8){$\scriptstyle 0$}
\put(16,8){$\scriptstyle  0$}
\put(0,-2){$\scriptstyle 0$}
\put(8,-2){$\scriptstyle 0$}
\put(14,-2){-$\scriptstyle  1$}
\dottedline{3}(0,5)(20,5)
\put(35,5){\vector(-1,0){10}}
\put(40,8){$\scriptstyle 0$}
\put(48,8){$\scriptstyle 0$}
\put(56,8){$\scriptstyle 0$}
\put(40,-2){$\scriptstyle 0$}
\put(46,-2){-$\scriptstyle 1$}
\put(56,-2){$\scriptstyle 0$}
\dottedline{3}(40,5)(60,5)
\put(65,5){\vector(1,0){10}}
\put(77, -5){\framebox(26,20)[c]}
\put(80,8){$\scriptstyle 0$}
\put(88,8){$\scriptstyle 0$}
\put(96,8){$\scriptstyle 0$}
\put(78,-2){-$\scriptstyle 1$}
\put(88,-2){$\scriptstyle 0$}
\put(96,-2){$\scriptstyle 0$}
\dottedline{3}(80,5)(100,5)%
\put(40,-15){$y(-4)$}
}
\put(150,270)
{
\put(-3, 35){\framebox(26,20)[c]}
\put(0,48){$\scriptstyle 0$}
\put(6,48){-$\scriptstyle 1$}
\put(14,48){-$\scriptstyle  1$}
\put(0,38){$\scriptstyle 0$}
\put(8,38){$\scriptstyle 0$}
\put(16,38){$\scriptstyle  0$}
\dottedline{3}(0,45)(20,45)%
\put(10,30){\vector(0,-1){10}}
\put(35,45){\vector(-1,0){10}}
\put(40,48){$\scriptstyle 0$}
\put(48,48){$\scriptstyle 1$}
\put(56,48){$\scriptstyle 0$}
\put(40,38){$\scriptstyle 0$}
\put(48,38){$\scriptstyle 0$}
\put(56,38){$\scriptstyle 0$}
\dottedline{3}(40,45)(60,45)
\put(50,20){\vector(0,1){10}}
\put(65,45){\vector(1,0){10}}
\put(77, 35){\framebox(26,20)[c]}
\put(78,48){-$\scriptstyle 1$}
\put(86,48){-$\scriptstyle 1$}
\put(96,48){$\scriptstyle 0$}
\put(80,38){$\scriptstyle 0$}
\put(88,38){$\scriptstyle 0$}
\put(96,38){$\scriptstyle 0$}
\dottedline{3}(80,45)(100,45)
\put(90,30){\vector(0,-1){10}}
%
\put(0,8){$\scriptstyle 0$}
\put(8,8){$\scriptstyle 0$}
\put(16,8){$\scriptstyle  0$}
\put(0,-2){$\scriptstyle 0$}
\put(8,-2){$\scriptstyle 0$}
\put(16,-2){$\scriptstyle  1$}
\dottedline{3}(0,5)(20,5)
\put(25,5){\vector(1,0){10}}
\put(37, -5){\framebox(26,20)[c]}
\put(40,8){$\scriptstyle 0$}
\put(48,8){$\scriptstyle 0$}
\put(56,8){$\scriptstyle 0$}
\put(38,-2){-$\scriptstyle 1$}
\put(46,-2){-$\scriptstyle 1$}
\put(54,-2){-$\scriptstyle 1$}
\dottedline{3}(40,5)(60,5)
\put(75,5){\vector(-1,0){10}}
%
\put(80,8){$\scriptstyle 0$}
\put(88,8){$\scriptstyle 0$}
\put(96,8){$\scriptstyle 0$}
\put(80,-2){$\scriptstyle 1$}
\put(88,-2){$\scriptstyle 0$}
\put(96,-2){$\scriptstyle 0$}
\dottedline{3}(80,5)(100,5)
\put(40,-15){$y(-3)$}
\put(-30,15){$\leftrightarrow$}
\put(-30,25){$\mu_+$}
}
\put(0,180)
{
%
\put(0,48){$\scriptstyle 0$}
\put(8,48){$\scriptstyle 1$}
\put(16,48){$\scriptstyle  1$}
\put(0,38){$\scriptstyle 0$}
\put(8,38){$\scriptstyle 0$}
\put(16,38){$\scriptstyle  0$}
\dottedline{3}(0,45)(20,45)%
%
\put(10,20){\vector(0,1){10}}
\put(25,45){\vector(1,0){10}}
\put(37, 35){\framebox(26,20)[c]}
\put(38,48){-$\scriptstyle 1$}
\put(46,48){-$\scriptstyle 1$}
\put(54,48){-$\scriptstyle 1$}
\put(40,38){$\scriptstyle 0$}
\put(48,38){$\scriptstyle 0$}
\put(56,38){$\scriptstyle 0$}
\dottedline{3}(40,45)(60,45)
\put(50,30){\vector(0,-1){10}}
\put(75,45){\vector(-1,0){10}}
%
\put(80,48){$\scriptstyle 1$}
\put(88,48){$\scriptstyle 1$}
\put(96,48){$\scriptstyle 0$}
\put(80,38){$\scriptstyle 0$}
\put(88,38){$\scriptstyle 0$}
\put(96,38){$\scriptstyle 0$}
\dottedline{3}(80,45)(100,45)
\put(90,20){\vector(0,1){10}}
\put(-3, -5){\framebox(26,20)[c]}
\put(0,8){$\scriptstyle 0$}
\put(8,8){$\scriptstyle 0$}
\put(16,8){$\scriptstyle  0$}
\put(-2,-2){-$\scriptstyle 1$}
\put(6,-2){-$\scriptstyle 1$}
\put(16,-2){$\scriptstyle  0$}
\dottedline{3}(0,5)(20,5)
\put(35,5){\vector(-1,0){10}}
\put(40,8){$\scriptstyle 0$}
\put(48,8){$\scriptstyle 0$}
\put(56,8){$\scriptstyle 0$}
\put(40,-2){$\scriptstyle 1$}
\put(48,-2){$\scriptstyle 1$}
\put(56,-2){$\scriptstyle 1$}
\dottedline{3}(40,5)(60,5)
\put(65,5){\vector(1,0){10}}
\put(77, -5){\framebox(26,20)[c]}
\put(80,8){$\scriptstyle 0$}
\put(88,8){$\scriptstyle 0$}
\put(96,8){$\scriptstyle 0$}
\put(80,-2){$\scriptstyle 0$}
\put(86,-2){-$\scriptstyle 1$}
\put(94,-2){-$\scriptstyle 1$}
\dottedline{3}(80,5)(100,5)%
\put(40,-15){$y(-2)$}
\put(-30,15){$\leftrightarrow$}
\put(-30,25){$\mu_-$}
}
\put(150,180)
{
\put(-3, 35){\framebox(26,20)[c]}
\put(-2,48){-$\scriptstyle 1$}
\put(8,48){$\scriptstyle 0$}
\put(16,48){$\scriptstyle  0$}
\put(0,38){$\scriptstyle 0$}
\put(8,38){$\scriptstyle 0$}
\put(16,38){$\scriptstyle  0$}
\dottedline{3}(0,45)(20,45)%
\put(10,30){\vector(0,-1){10}}
\put(35,45){\vector(-1,0){10}}
\put(40,48){$\scriptstyle 1$}
\put(48,48){$\scriptstyle 1$}
\put(56,48){$\scriptstyle 1$}
\put(40,38){$\scriptstyle 0$}
\put(48,38){$\scriptstyle 0$}
\put(56,38){$\scriptstyle 0$}
\dottedline{3}(40,45)(60,45)
\put(50,20){\vector(0,1){10}}
\put(65,45){\vector(1,0){10}}
\put(77, 35){\framebox(26,20)[c]}
\put(80,48){$\scriptstyle 0$}
\put(88,48){$\scriptstyle 0$}
\put(94,48){-$\scriptstyle 1$}
\put(80,38){$\scriptstyle 0$}
\put(88,38){$\scriptstyle 0$}
\put(96,38){$\scriptstyle 0$}
\dottedline{3}(80,45)(100,45)
\put(90,30){\vector(0,-1){10}}
%
\put(0,8){$\scriptstyle 0$}
\put(8,8){$\scriptstyle 0$}
\put(16,8){$\scriptstyle  0$}
\put(0,-2){$\scriptstyle 1$}
\put(8,-2){$\scriptstyle 1$}
\put(16,-2){$\scriptstyle  0$}
\dottedline{3}(0,5)(20,5)
\put(25,5){\vector(1,0){10}}
\put(37, -5){\framebox(26,20)[c]}
\put(40,8){$\scriptstyle 0$}
\put(48,8){$\scriptstyle 0$}
\put(56,8){$\scriptstyle 0$}
\put(40,-2){$\scriptstyle 0$}
\put(46,-2){-$\scriptstyle 1$}
\put(56,-2){$\scriptstyle 0$}
\dottedline{3}(40,5)(60,5)
\put(75,5){\vector(-1,0){10}}
%
\put(80,8){$\scriptstyle 0$}
\put(88,8){$\scriptstyle 0$}
\put(96,8){$\scriptstyle 0$}
\put(80,-2){$\scriptstyle 0$}
\put(88,-2){$\scriptstyle 1$}
\put(96,-2){$\scriptstyle 1$}
\dottedline{3}(80,5)(100,5)
\put(40,-15){$y(-1)$}
\put(-30,15){$\leftrightarrow$}
\put(-30,25){$\mu_+$}
}
\put(0,90)
{
%
\put(0,48){$\scriptstyle 1$}
\put(8,48){$\scriptstyle 0$}
\put(16,48){$\scriptstyle  0$}
\put(0,38){$\scriptstyle 0$}
\put(8,38){$\scriptstyle 0$}
\put(16,38){$\scriptstyle  0$}
%
\put(10,20){\vector(0,1){10}}
\put(25,45){\vector(1,0){10}}
\put(37, 35){\framebox(26,20)[c]}
\put(40,48){$\scriptstyle 0$}
\put(48,48){$\scriptstyle 1$}
\put(56,48){$\scriptstyle 0$}
\put(40,38){$\scriptstyle 0$}
\put(48,38){$\scriptstyle 0$}
\put(56,38){$\scriptstyle 0$}
%
\put(50,30){\vector(0,-1){10}}
\put(75,45){\vector(-1,0){10}}
%
\put(80,48){$\scriptstyle 0$}
\put(88,48){$\scriptstyle 0$}
\put(96,48){$\scriptstyle 1$}
\put(80,38){$\scriptstyle 0$}
\put(88,38){$\scriptstyle 0$}
\put(96,38){$\scriptstyle 0$}
%
\put(90,20){\vector(0,1){10}}
\put(-3, -5){\framebox(26,20)[c]}
\put(0,8){$\scriptstyle 0$}
\put(8,8){$\scriptstyle 0$}
\put(16,8){$\scriptstyle  0$}
\put(0,-2){$\scriptstyle 1$}
\put(8,-2){$\scriptstyle 0$}
\put(16,-2){$\scriptstyle  0$}
%
\put(35,5){\vector(-1,0){10}}
\put(40,8){$\scriptstyle 0$}
\put(48,8){$\scriptstyle 0$}
\put(56,8){$\scriptstyle 0$}
\put(40,-2){$\scriptstyle 0$}
\put(48,-2){$\scriptstyle 1$}
\put(56,-2){$\scriptstyle 0$}
%
\put(65,5){\vector(1,0){10}}
\put(77, -5){\framebox(26,20)[c]}
\put(80,8){$\scriptstyle 0$}
\put(88,8){$\scriptstyle 0$}
\put(96,8){$\scriptstyle 0$}
\put(80,-2){$\scriptstyle 0$}
\put(88,-2){$\scriptstyle 0$}
\put(96,-2){$\scriptstyle 1$}
%
\put(40,-15){$y(0)$}
\put(-30,15){$\leftrightarrow$}
\put(-30,25){$\mu_-$}
}
\put(150,90)
{
\put(-3, 35){\framebox(26,20)[c]}
\put(0,48){$\scriptstyle 1$}
\put(8,48){$\scriptstyle 0$}
\put(16,48){$\scriptstyle  0$}
\put(0,38){$\scriptstyle 1$}
\put(8,38){$\scriptstyle 0$}
\put(16,38){$\scriptstyle  0$}
\dottedline{3}(5,35)(5,55)
\dottedline{3}(13,35)(13,55)
%
\put(10,30){\vector(0,-1){10}}
\put(35,45){\vector(-1,0){10}}
\put(40,48){$\scriptstyle 0$}
\put(46,48){-$\scriptstyle 1$}
\put(56,48){$\scriptstyle 0$}
\put(40,38){$\scriptstyle 0$}
\put(48,38){$\scriptstyle 0$}
\put(56,38){$\scriptstyle 0$}
\dottedline{3}(45,35)(45,55)
\dottedline{3}(53,35)(53,55)
%
\put(50,20){\vector(0,1){10}}
\put(65,45){\vector(1,0){10}}
\put(77, 35){\framebox(26,20)[c]}
\put(80,48){$\scriptstyle 0$}
\put(88,48){$\scriptstyle 0$}
\put(96,48){$\scriptstyle 1$}
\put(80,38){$\scriptstyle 0$}
\put(88,38){$\scriptstyle 0$}
\put(96,38){$\scriptstyle 1$}
\dottedline{3}(85,35)(85,55)
\dottedline{3}(93,35)(93,55)
%
\put(90,30){\vector(0,-1){10}}
%
\put(0,8){$\scriptstyle 0$}
\put(8,8){$\scriptstyle 0$}
\put(16,8){$\scriptstyle  0$}
\put(-2,-2){-$\scriptstyle 1$}
\put(8,-2){$\scriptstyle 0$}
\put(16,-2){$\scriptstyle  0$}
\dottedline{3}(5,-5)(5,15)
\dottedline{3}(13,-5)(13,15)
%
\put(25,5){\vector(1,0){10}}
\put(37, -5){\framebox(26,20)[c]}
\put(40,8){$\scriptstyle 0$}
\put(48,8){$\scriptstyle 1$}
\put(56,8){$\scriptstyle 0$}
\put(40,-2){$\scriptstyle 0$}
\put(48,-2){$\scriptstyle 1$}
\put(56,-2){$\scriptstyle 0$}
\dottedline{3}(45,-5)(45,15)
\dottedline{3}(53,-5)(53,15)
%
\put(75,5){\vector(-1,0){10}}
%
\put(80,8){$\scriptstyle 0$}
\put(88,8){$\scriptstyle 0$}
\put(96,8){$\scriptstyle 0$}
\put(80,-2){$\scriptstyle 0$}
\put(88,-2){$\scriptstyle 0$}
\put(94,-2){-$\scriptstyle 1$}
\dottedline{3}(85,-5)(85,15)
\dottedline{3}(93,-5)(93,15)
%
\put(40,-15){$y(1)$}
\put(-30,15){$\leftrightarrow$}
\put(-30,25){$\mu_+$}
}
\put(0,0)
{
%
\put(-2,48){-$\scriptstyle 1$}
\put(8,48){$\scriptstyle 0$}
\put(16,48){$\scriptstyle  0$}
\put(-2,38){-$\scriptstyle 1$}
\put(8,38){$\scriptstyle 0$}
\put(16,38){$\scriptstyle  0$}
\dottedline{3}(5,35)(5,55)
\dottedline{3}(13,35)(13,55)
%
\put(10,20){\vector(0,1){10}}
\put(25,45){\vector(1,0){10}}
\put(37, 35){\framebox(26,20)[c]}
\put(40,48){$\scriptstyle 0$}
\put(48,48){$\scriptstyle 0$}
\put(56,48){$\scriptstyle 0$}
\put(40,38){$\scriptstyle 0$}
\put(48,38){$\scriptstyle 1$}
\put(56,38){$\scriptstyle 0$}
\dottedline{3}(45,35)(45,55)
\dottedline{3}(53,35)(53,55)
\put(50,30){\vector(0,-1){10}}
\put(75,45){\vector(-1,0){10}}
%
\put(80,48){$\scriptstyle 0$}
\put(88,48){$\scriptstyle 0$}
\put(94,48){-$\scriptstyle 1$}
\put(80,38){$\scriptstyle 0$}
\put(88,38){$\scriptstyle 0$}
\put(94,38){-$\scriptstyle 1$}
\dottedline{3}(85,35)(85,55)
\dottedline{3}(93,35)(93,55)
\put(90,20){\vector(0,1){10}}
\put(-3, -5){\framebox(26,20)[c]}
\put(0,8){$\scriptstyle 1$}
\put(8,8){$\scriptstyle 0$}
\put(16,8){$\scriptstyle  0$}
\put(0,-2){$\scriptstyle 0$}
\put(8,-2){$\scriptstyle 0$}
\put(16,-2){$\scriptstyle  0$}
\dottedline{3}(5,-5)(5,15)
\dottedline{3}(13,-5)(13,15)
\put(35,5){\vector(-1,0){10}}
\put(40,8){$\scriptstyle 0$}
\put(46,8){-$\scriptstyle 1$}
\put(56,8){$\scriptstyle 0$}
\put(40,-2){$\scriptstyle 0$}
\put(46,-2){-$\scriptstyle 1$}
\put(56,-2){$\scriptstyle 0$}
\dottedline{3}(45,-5)(45,15)
\dottedline{3}(53,-5)(53,15)
\put(65,5){\vector(1,0){10}}
\put(77, -5){\framebox(26,20)[c]}
\put(80,8){$\scriptstyle 0$}
\put(88,8){$\scriptstyle 0$}
\put(96,8){$\scriptstyle 1$}
\put(80,-2){$\scriptstyle 0$}
\put(88,-2){$\scriptstyle 0$}
\put(96,-2){$\scriptstyle 0$}
\dottedline{3}(85,-5)(85,15)
\dottedline{3}(93,-5)(93,15)
\put(40,-15){$y(2)$}
\put(-30,15){$\leftrightarrow$}
\put(-30,25){$\mu_-$}
}
\put(150,0)
{
\put(-3, 35){\framebox(26,20)[c]}
\put(0,48){$\scriptstyle 0$}
\put(8,48){$\scriptstyle 0$}
\put(16,48){$\scriptstyle  0$}
\put(-2,38){-$\scriptstyle 1$}
\put(8,38){$\scriptstyle 0$}
\put(16,38){$\scriptstyle  0$}
\dottedline{3}(5,35)(5,55)
\dottedline{3}(13,35)(13,55)
\put(10,30){\vector(0,-1){10}}
\put(35,45){\vector(-1,0){10}}
\put(40,48){$\scriptstyle 0$}
\put(48,48){$\scriptstyle 0$}
\put(56,48){$\scriptstyle 0$}
\put(40,38){$\scriptstyle 0$}
\put(46,38){-$\scriptstyle 1$}
\put(56,38){$\scriptstyle 0$}
\dottedline{3}(45,35)(45,55)
\dottedline{3}(53,35)(53,55)
\put(50,20){\vector(0,1){10}}
\put(65,45){\vector(1,0){10}}
\put(77, 35){\framebox(26,20)[c]}
\put(80,48){$\scriptstyle 0$}
\put(88,48){$\scriptstyle 0$}
\put(96,48){$\scriptstyle 0$}
\put(80,38){$\scriptstyle 0$}
\put(88,38){$\scriptstyle 0$}
\put(94,38){-$\scriptstyle 1$}
\dottedline{3}(85,35)(85,55)
\dottedline{3}(93,35)(93,55)
\put(90,30){\vector(0,-1){10}}
%
\put(-2,8){-$\scriptstyle 1$}
\put(8,8){$\scriptstyle 0$}
\put(16,8){$\scriptstyle  0$}
\put(0,-2){$\scriptstyle 0$}
\put(8,-2){$\scriptstyle 0$}
\put(16,-2){$\scriptstyle  0$}
\dottedline{3}(5,-5)(5,15)
\dottedline{3}(13,-5)(13,15)
\put(25,5){\vector(1,0){10}}
\put(37, -5){\framebox(26,20)[c]}
\put(40,8){$\scriptstyle 0$}
\put(46,8){-$\scriptstyle 1$}
\put(56,8){$\scriptstyle 0$}
\put(40,-2){$\scriptstyle 0$}
\put(48,-2){$\scriptstyle 0$}
\put(56,-2){$\scriptstyle 0$}
\dottedline{3}(45,-5)(45,15)
\dottedline{3}(53,-5)(53,15)
\put(75,5){\vector(-1,0){10}}
%
\put(80,8){$\scriptstyle 0$}
\put(88,8){$\scriptstyle 0$}
\put(94,8){-$\scriptstyle 1$}
\put(80,-2){$\scriptstyle 0$}
\put(88,-2){$\scriptstyle 0$}
\put(96,-2){$\scriptstyle 0$}
\dottedline{3}(85,-5)(85,15)
\dottedline{3}(93,-5)(93,15)
\put(40,-15){$y(3)$}
\put(-30,15){$\leftrightarrow$}
\put(-30,25){$\mu_+$}
}
\end{picture}
\end{center}

A closer look at the diagram reveals the following
{\em factorization property\/}  \cite{Nakanishi09}.

\begin{itemize}
\item[(i)] In the region $u=0,1,2$, the variables $[y_i(u)]$
at the forward mutation points (framed one) transform
like the positive roots of $A_2$ in each column.
Compare it with the second diagram in
 Example \ref{ex:a2l2}.
This is because during this region they remain positive
and then, by the sign-arrow coordination \eqref{eq:signarrow},
the {\em horizontal\/} arrows can be ignored in the mutation.

\item[(ii)]  In the region $u=-1,-2,-3,-4$, the variables $[y_i(u)]$
at the forward mutation points transform
like the positive roots of $A_3$  in each column.
Compare it with two diagrams in
 Example \ref{ex:a3l2}.
Again, this is because during this region they remain negative
and then, by the sign-arrow coordination \eqref{eq:signarrow},
the {\em vertical\/} arrows can be ignored in the mutation.
\end{itemize}
This explains the formula of periodicity $7=4+3$.
One can easily generalize the argument to any
 simply laced $X$ and any $\ell$, using the result of \cite{Fomin03b}
for $\ell=2$.
\end{example}

The periodicities of seeds corresponding to the Y-systems of nonsimply laced
type can be proved in the same spirit, but they are necessarily
more complicated
(thus, perhaps more intriguing).
See \cite{Inoue10a,Inoue10b}.

\section{Application II: Dilogarithm identities}

\subsection{Dilogarithm functions}

Define the Euler dilogarithm
$\mathrm{Li}_2(x)$ and Rogers dilogarithm $L(x)$
by the following integrals \cite{Lewin81,Kirillov95,Zagier07}.
\begin{align}
\label{eq:L00}
\mathrm{Li}_2(x)
&=-\int_{0}^x
\left\{ \frac{\log(1-y)}{y}
\right\} dy
\quad (x\leq 1),
\\
\label{eq:L0}
L(x)&=-\frac{1}{2}\int_{0}^x
\left\{ \frac{\log(1-y)}{y}+
\frac{\log y}{1-y}
\right\} dy
\quad (0 \leq x\leq 1).
\end{align}
We restrict the argument $x$ in the above to
avoid the branches coming from the logarithm.
Two functions are related by
\begin{align}
\label{eq:LL1}
L(x)
&=\mathrm{Li}_2(x)
+ \frac{1}{2}
\log x \log (1-x)
\quad (0 \leq x \leq 1),
\\
\label{eq:LL2}
-L\left(\frac{x}{1+x}\right)
&=\mathrm{Li}_2(-x)
+ \frac{1}{2}
 \log x \log (1+x)
\quad (0 \leq x).
\end{align}
The Rogers dilogarithm satisfies the following properties:
\begin{align}
L(0)&=0,\quad L(1)=\zeta(2)=\frac{\pi^2}{6},\\
\label{eq:euler}
L(x)&+L(1-x)=\frac{\pi^2}{6}
\quad (0 \leq x \leq 1),\\
\label{eq:Lpent}
\begin{split}
L(x)&+L(y)
+L\left(\frac{1-x}{1-xy}\right)
+
L(1-xy)+
L\left(\frac{1-y}{1-xy}\right)
=\frac{\pi^2}{2}\\
&\hskip160pt (0\leq x,y\leq 1).
\end{split}
\end{align}
The relations \eqref{eq:euler} and \eqref{eq:Lpent}
are called Euler's identity and Abel's identity
(or the pentagon identity).
Using \eqref{eq:euler}, the pentagon identity is also written 
as follows:
\begin{align}
\label{eq:Lpent2}
L(x)+L(y)=
L\left(\frac{x(1-y)}{1-xy}\right)
+
L(xy)+
L\left(\frac{y(1-x)}{1-xy}\right).
\end{align}

\subsection{Background: Dilogarithm identities in CFT}

In the late 80's there appeared a remarkable
conjecture on the dilogarithm identities
of central charges in conformal field theory (CFT)
\cite{Kirillov89,Kirillov90,Bazhanov90,Kuniba93a}.
They  are associated with
pairs $(X,\ell)$ of a Dynkin diagram $X$ of
finite type and an integer $\ell\geq 2$  (level),
which are the same data for the Y-systems
in Section \ref{subsec:period}.

For simplicity, let us concentrate
on the case when $X$ is simply laced case
as in Section \ref{subsec:period}.
For a family of variables
$\{ Y^{(a)}_m \mid
a\in \mathcal{I};
m=1,\dots,\ell-1 \}$,
we introduce the following system of
algebraic equations, which is called the {\em constant
Y-system}:
\begin{align}
\label{eq:y1c}
(Y^{(a)}_m)^2
&=
\frac{
\prod_{b\in I, b\sim a}
(1+Y^{(b)}_m)
}
{
(1+Y^{(a)}_{m-1}{}^{-1})(1+Y^{(a)}_{m+1}{}^{-1})
},
\end{align}
where $Y^{(a)}_{0}{}^{-1}=Y^{(a)}_{\ell}{}^{-1}=0$.
Namely, it is the Y-system \eqref{eq:Ysystem}
with the condition that the variables $Y^{(a)}_m(u)$
are constant with respect to $u$.

\begin{theorem}[\cite{Nahm09}]
There exists a unique real positive solution of \eqref{eq:y1c}.
\end{theorem}

The dilogarithm conjecture claimed that
for the above real positive solution of \eqref{eq:y1c},
the following equality holds.
\begin{align}
\label{eq:DI6}
 \frac{6}{\pi^2}
\sum_{a\in \mathcal{I}}
\sum_{m=1}^{\ell-1}
L\left(
\frac{Y^{(a)}_m}
{1+Y^{(a)}_m}
\right)=
\frac{\ell \dim \mathfrak{g}}{h+\ell}-n,
\end{align}
where $h$ is the Coxeter number of $X$,
 $n=|\mathcal{I}|= \mathrm{rank}\, X$, 
 and $\mathfrak{g}$ is the simple Lie algebra of type
 $X$.
The equality \eqref{eq:DI6} was proved 
for $X=A_n$
by \cite{Kirillov89} by using an explicit solution
of \eqref{eq:y1c}.
The first term of the right hand side
of \eqref{eq:DI6} is the central charge of
the {\em Wess-Zumino-Witten model} of type $X$ and level $\ell$.
The right hand side itself is
the central charge of the {\em parafermion CFT}
of type $X$ and level $\ell$.
See \cite{Kuniba10} and references therein for further background
of the equality \eqref{eq:DI6}.

Gliozzi and Tateo proposed the functional
generalization of 
\eqref{eq:DI6}
in accordance with the periodicity conjecture
\eqref{eq:yperiod} of the Y-systems   \cite{Gliozzi95}.
Namely, for any real positive solution of
the Y-system \eqref{eq:Ysystem},
the following identity holds.
\begin{align}
\label{eq:DI7}
 \frac{6}{\pi^2}
\sum_{a\in \mathcal{I}}
\sum_{m=1}^{\ell-1}
\sum_{u=0}^{2(h + \ell)-1}
L\left(
\frac{Y^{(a)}_m(u)}
{1+Y^{(a)}_m(u)}
\right)=2hn(\ell-1).
\end{align}
Indeed, the equality \eqref{eq:DI6} follows from 
\eqref{eq:DI7} by specializing it to the constant
solution with respect to $u$.
The identity \eqref{eq:DI7} was proved for $X=A_1$ by
\cite{Gliozzi96} and 
\cite{Frenkel95} by using  explicit solutions of
the Y-systems.
It was proved for simply laced
$X$ and $\ell=2$ by \cite{Chapoton05} using the cluster algebra
method based on the result of \cite{Fomin03a}.
Then, it was proved in full generality by
\cite{Nakanishi09,Inoue10a,Inoue10b} using the tropicalization
method.
It was further generalized to the dilogarithm identities
associated with periods of seeds in cluster algebras
\cite{Nakanishi10c}.

\subsection{Classical dilogarithm identities}
{\em In this subsection we assume that $B$ is skew-symmetric.}

Let us present so far the most general dilogarithm identities
associated with periods of seeds in cluster algebras
\cite{Nakanishi10c}.
Since $x$-variables are irrelevant, let us concentrate
on `$y$-seeds' $(B',y')$.
Let $(k_0,\dots,k_{N-1})$ be a $\nu$-period of the initial
seed $(B,y)$.
Set$(B(0),y(0))=(B,y)$,
and consider the sequence of mutations,
\begin{align}
\label{eq:seedmutseq}
&(B(0),y(0))
\
\mathop{\leftrightarrow}^{\mu_{k_0}}
\
(B(1),y(1))
\
\mathop{\leftrightarrow}^{\mu_{k_1}}
\
\cdots
\
\mathop{\leftrightarrow}^{\mu_{k_{N-1}}}
\
(B(N),y(N))
\end{align}
Let $\varepsilon_t$ be the tropical sign 
of $(B(t),y(t))$ at $k_t$.
We call the sequence
 $(\varepsilon_0,
\dots, \varepsilon_{N-1})$
the {\em tropical sign-sequence}
of \eqref{eq:seedmutseq}.

\begin{theorem}[Classical dilogarithm identities \cite{Nakanishi10c}]
\label{thm:CDI2}
The following equalities hold for any evaluation
of the initial $y$-variables $y_i$ $(i\in I)$ in $\mathbb{R}_{>0}$.
\begin{align}
\label{eq:DI3}
\sum_{t=0}^{N-1}
\varepsilon_t
L\left(
\frac{y_{k_t}(t)^{\varepsilon_t}}
{1+y_{k_t}(t)^{\varepsilon_t}}
\right)&=0,\\
\label{eq:DI4}
 \frac{6}{\pi^2}
\sum_{t=0}^{N-1}
L\left(
\frac{y_{k_t}(t)}
{1+y_{k_t}(t)}
\right)&=N_-,\\
\label{eq:DI10}
 \frac{6}{\pi^2}
\sum_{t=0}^{N-1}
L\left(
\frac{1}
{1+y_{k_t}(t)}
\right)&=N_+,
\end{align}
where $N_+$ and $N_-$ ($N_+ + N_- = N$) are the total numbers
of
 $1$ and $-1$ among $\varepsilon_0,\dots,\varepsilon_{N-1}$,
respectively. 
\end{theorem}

Three identities are equivalent to each other
due to \eqref{eq:euler}.
Applying \eqref{eq:DI4} for the periods corresponding to
the Y-systems 
and counting the number $N_-$, we obtain
\eqref{eq:DI7}, thus proving the dilogarithm identities
\eqref{eq:DI6} as well.

\begin{example}[Pentagon identity]
\label{ex:pent}
Apply Theorem \ref{thm:CDI2} to the period in
Example \ref{ex:a2}.
We have $(k_0,k_1,k_2,k_3,k_4)=(1,2,1,2,1)$ and
\begin{align}
\begin{split}
&y_1(0)=y_1,
\quad
y_2(1)=y_2(1+y_1),
\quad
y_1(2)=y_1^{-1}(1+y_2+y_1y_2),\\
&
y_2(3)=y_1^{-1}y_2^{-1}(1+y_2),
\quad
y_1(4)=y_2^{-1},\\
&\varepsilon_0 =
\varepsilon_1 = 1,
\quad
\varepsilon_2 = 
\varepsilon_3 = \varepsilon_4 = -1.
\end{split}
\end{align}
Inserting these data into \eqref{eq:DI3}, we obtain
\begin{align}
\begin{split}
\label{eq:DI5}
L&\left(
\frac{y_{1}}
{1+y_{1}}
\right)
+
L\left(
\frac{y_{2}(1+y_1)}
{1+y_{2}+y_1y_2}
\right)
\\
&-
L\left(
\frac{y_1}
{(1+y_1)(1+y_2)}
\right)
-
L\left(
\frac{y_1y_2}
{1+y_{2}+y_1y_2}
\right)
-
L\left(
\frac{y_{2}}
{1+y_{2}}
\right)
=0.
\end{split}
\end{align}
By setting
$x=y_1/(1+y_1)$, $y=y_2(1+y_2)/(1+y_2+y_1y_2)$,
it coincides with the pentagon identity \eqref{eq:Lpent2}.
\end{example}

\subsection{Quantum pentagon identity}

Following
\cite{Faddeev93,Faddeev94},
define the {\em quantum dilogarithm $\mathbf{\Psi}_q(x)$},
 for $|q|<1$ and $x\in \mathbb{C}$, by
\begin{align}
\label{eq:psi1}
\mathbf{\Psi}_q(x)&=
\sum_{n=0}^\infty
\frac{(-qx)^n}{(q^2;q^2)_n}=
\frac{1}{(-qx;q^2)_{\infty}},
\quad
(a;q)_n = \prod_{k=0}^{n-1}(1-aq^{k}).
\end{align}
It satisfies the following recursion relations.
\begin{align}
\label{eq:rec1}
\mathbf{\Psi}_q(q^{\pm2}x)& =
(1+q^{\pm 1}x)^{\pm1}\mathbf{\Psi}_q(x).
\end{align}
The following properties explain why
 it is considered as a quantum analogue of  the dilogarithm  \cite{Faddeev93,Faddeev94}.

{\em
$(a)$. Asymptotic behavior:
  In  the limit $q \rightarrow 1^{-}$,
\begin{align}
\label{eq:asym3}
\mathbf{\Psi}_q(x) \sim
\exp \left( -\frac{\mathrm{Li}_2(-x)}{2 \log q }
\right).
\end{align}

$(b)$. Quantum pentagon identity: If $UV=q^2VU$, then
\begin{align}
\label{eq:pent2}
\mathbf{\Psi}_q(U)\mathbf{\Psi}_q(V)
=
\mathbf{\Psi}_q(V)\mathbf{\Psi}_q(q^{-1}UV)\mathbf{\Psi}_q(U).
\end{align}
Moreover, in the limit $q\rightarrow 1^{-}$, the relation
\eqref{eq:pent2} reduces to the   pentagon identity
\eqref{eq:Lpent2}.
}

\subsection{Quantum dilogarithm identities}
{\em In this subsection we assume that $B$ is skew-symmetric.}

We present a quantum counterpart of
the dilogarithm identities \eqref{eq:DI3},
which is a generalization of \eqref{eq:pent2}.
To do that, we use
the quantum cluster algebras by \cite{Fock03,Fock07}.
In short, it replaces classical $y$-seeds
$(B',y')$ with the quantum one $(B',Y')$ in the following way:
First, the initial quantum $y$-variables are
 {\em noncommutative\/} and satisfy the relation
\begin{align}
\label{eq:yy1}
Y_i Y_j = q^{2b_{ji}} Y_j Y_i.
\end{align}
Second, for the mutation $(B'',Y'')=\mu_k(B',Y')$,
the $\varepsilon$-expression of the exchange relation
\eqref{eq:ymut2} is replaced with
\begin{align}
\label{eq:coefq}
Y''_i =
\begin{cases}
\displaystyle
{Y'_k}{}^{-1}&i=k\\
\displaystyle
q^{
b'_{ik}[\varepsilon b'_{ki}]_+
}
Y'_i Y'_k{}^{[\varepsilon b'_{ki}]_+}
\prod_{m=1}^{|b'_{ki}|}
(1+ q^{-\varepsilon \mathrm{sgn}(b'_{ki})(2m-1)}
Y'_k{}^{\varepsilon})^{-\mathrm{sgn}({b'_{ki}})}
&
i\neq k.\\
\end{cases}
\end{align}
Formally setting $q=1$, quantum $y$-seeds reduce
to the classical one.

Let  $\varepsilon'_k$ be the tropical sign  of $(B',y')$ at $k$.
In analogy with the classical case,
 we introduce the {\em tropical quantum $y$-variables
$[Y'_i]$} by the initial condition $[Y_i]=Y_i$ and the exchange relation
\begin{align}
\label{eq:tropcoefq}
[Y''_i] =
\begin{cases}
\displaystyle
[{Y'_k}]{}^{-1}&i=k\\
\displaystyle
q^{
b'_{ik}[\varepsilon'_k b'_{ki}]_+
}
[Y'_i][ Y'_k]{}^{[\varepsilon'_k b'_{ki}]_+}
&
i\neq k.\\
\end{cases}
\end{align}
This gives the tropical part of 
 \eqref{eq:coefq} with $\varepsilon=\varepsilon'_k$ therein.
 
 On the other hand, the nontropical part of  \eqref{eq:coefq} 
 is given by the adjoint action of the quantum dilogarithm;
 More precisely,
\begin{align}
\label{eq:ad}
\begin{split}
\mathsf{Ad}(\mathbf{\Psi}_q(\mathsf{Y}'_k{}^{\varepsilon'_k}){}^{\varepsilon'_k})
(\mathsf{Y}'_i) := &
\ \mathbf{\Psi}_q(\mathsf{Y}'_k{}^{\varepsilon'_k}) {}^{\varepsilon'_k}
\mathsf{Y}'_i
\mathbf{\Psi}_q(\mathsf{Y}'_k{}^{\varepsilon'_k}){}^{-\varepsilon'_k}\\
=& \ \mathsf{Y}'_i
\mathbf{\Psi}_q(q^{-2b'_{ki}}\mathsf{Y}'_k{}^{\varepsilon'_k}){}^{\varepsilon'_k}
\mathbf{\Psi}_q(\mathsf{Y}'_k{}^{\varepsilon'_k}){}^{-\varepsilon'_k}\\
=&\
\mathsf{Y}'_i \prod_{m=1}^{|b'_{ki}|}
(1+ q^{-\varepsilon'_k\mathrm{sgn}(b'_{ki})(2m-1)}
\mathsf{Y}'_k{}^{\varepsilon'_k})^{-\mathrm{sgn}({b'_{ki}})},
\end{split}
\end{align}
where we use \eqref{eq:rec1} in the last equality.
 This is where the quantum dilogarithm is involved in quantum cluster algebras.
Actually, in \cite{Fock03} this important property 
(without $\varepsilon'_k$) was
employed as the definition of the exchange relation of the
quantum $y$-variables.
The factor $\varepsilon'_k$ was introduced by
 \cite{Keller11} based on the work \cite{Nagao10}.

A $\nu$-period of a quantum $y$-seed is defined in the same way as
the classical case.
One can show that
an $I$-sequence  is a $\nu$-period of a quantum $y$-seed $(B',Y')$ 
if and only if it is a $\nu$-period of the corresponding $y$-seed $(B',y')$
\cite[Prop.3.4]{Kashaev11}, using the result of
 \cite[Theorem 6.1]{Berenstein05b} and \cite[Lemma 2.22]{Fock07}.

\begin{theorem}[Quantum dilogarithm identities \cite{Keller11,Nagao11b}]
\label{thm:QDI1}
Suppose that  $(k_0,\dots,k_{N-1})$ is a $\nu$-period
of $(B,Y)$,
and let  $(\varepsilon_1$,
\dots, $ \varepsilon_{N-1})$  be the tropical sign-sequence
of \eqref{eq:seedmutseq}.
 Then, the following identity holds.
\begin{align}
\label{eq:QDI2}
\mathbf{\Psi}_q([Y_{k_0}(0)]^{\varepsilon_0}
 )^{\varepsilon_0}
\cdots
\mathbf{\Psi}_q([Y_{k_{N-1}}(N-1)]^{\varepsilon_{N-1} }
 )^{\varepsilon_{N-1}}
=1.
\end{align}
\end{theorem}

\begin{example}[Quantum pentagon identity]
\label{ex:qpent}
We continue to use the data in
Examples \ref{ex:a2} and \ref{ex:pent}.
For the initial quantum $y$-seed $(B,Y)$, we have
\begin{align}
Y_1 Y_2 = q^{2} Y_2 Y_1.
\end{align}
For  the quantum $y$-seeds
corresponding to \eqref{eq:seedmutseq3},
we have
\begin{align}
\begin{split}
\begin{cases}
Y_1(0)=Y_1\\
Y_2(0)=Y_2,\\
\end{cases}
\begin{cases}
Y_1(1)=Y_1^{-1}\\
Y_2(1)=Y_2(1+q Y_1),\\
\end{cases}
\begin{cases}
Y_1(2)=Y_1^{-1}(1+q Y_2 + Y_1Y_2)\\
Y_2(2)=Y_2^{-1}(1+q^{-1} Y_1)^{-1},\\
\end{cases}\\
\begin{cases}
Y_1(3)=Y_1(1+q Y_2 + Y_1Y_2)^{-1}\\
Y_2(3)=q ^{-1}Y_1^{-1}Y_2^{-1} (1+q Y_2),\\
\end{cases}
\begin{cases}
Y_1(4)=Y_2^{-1}\\
Y_2(4)=q^{-1}Y_1Y_2(1+q^{-1} Y_2),\\
\end{cases}
\begin{cases}
Y_1(5)=Y_2\\
Y_2(5)=Y_1,\\
\end{cases}
\end{split}
\end{align}
\begin{align}
\begin{split}
[Y_1(0)]&=Y_1,
\quad
[Y_2(1)]=Y_2,
\quad
[Y_1(2)]=Y_1^{-1},\\
\quad
[Y_2(3)]&=q^{-1}Y_1^{-1}Y_2^{-1}=(qY_2 Y_1)^{-1}=(q^{-1}Y_1 Y_2)^{-1},
\quad
[Y_1(4)]=Y_2^{-1}.
\end{split}
\end{align}
Thus, the quantum dilogarithm identity  \eqref{eq:QDI2}
is
\begin{align}
\begin{split}
\label{eq:QDI6}
\mathbf{\Psi}_q\left(
Y_1
\right)
\mathbf{\Psi}_q\left(
Y_2
\right)
\mathbf{\Psi}_q\left(
Y_1
\right)^{-1}
\mathbf{\Psi}_q\left(
q^{-1} Y_1 Y_2
\right)^{-1}
\mathbf{\Psi}_q\left(
Y_2
\right)^{-1}
=1
\end{split}.
\end{align}
It coincides with
the quantum pentagon relation \eqref{eq:pent2}.
\end{example}

In \cite{Kashaev11} several forms of quantum dilogarithm identities
are given,
where the identities \eqref{eq:QDI2} are called the {\em tropical form}.
Rewriting them in the {\em local form} therein,
then taking the semiclassical limit $q\rightarrow 1$
with the saddle point method, 
the classical dilogarithm identities  \eqref{eq:DI3}
can be recovered \cite{Kashaev11}.

In summary, we obtain the following scheme,
which vastly generalizes the classical and quantum pentagon identities:
\begin{align*}
\begin{xy}
\xymatrix{
\framebox{periods of quantum cluster algebras}
\ar@{->}[r]&
\framebox{quantum dilogarithm identities}
\ar@{->}[d]^{\mbox{semiclassical limit}}
\\
\framebox{periods of classical cluster algebras}
\ar@{->}[r] \ar@{<->}[u] &
\framebox{classical dilogarithm identities}
\\
}
\end{xy}
\end{align*}

\bibliography{../../biblist/biblist11-1.bib}

\end{document}